%% file: ex_article.tex
\documentclass[hidelinks,onefignum,onetabnum]{siamart251216}

\input{ex_shared}

\ifpdf
\hypersetup{
	pdftitle={Sparsely connected rank-inspired neural network},
	pdfauthor={Y. Huang, W. Su, N. Yi and P. Yin}
}
\fi

\begin{document}
	
	\maketitle
	
	% REQUIRED
	\begin{abstract}
		We propose a Sparsely Connected Rank-Inspired Neural Network (SCRINN) for function approximation and the numerical solution of partial differential equations (PDEs). By introducing a structured sparse connectivity architecture and an orthogonality- and normalization-guided loss function, SCRINN constructs an approximately orthonormal neural basis while promoting progressive rank growth of its evaluation matrix and reducing the number of trainable parameters. A two-stage framework is developed, in which an approximately orthonormal neural basis is first constructed, and the output weights are subsequently determined using least-squares methods. Numerical experiments on function approximation and a variety of steady-state, time-dependent, and nonlinear PDEs demonstrate that SCRINN consistently achieves high approximation accuracy with substantially lower training costs. The code is publicly available on GitHub.\footnotemark
	\end{abstract}
	\footnotetext{GitHub: \url{https://github.com/suwu-suwu/SCRINN}.}
	% REQUIRED
	
	% REQUIRED

	\section{Introduction}
	Function approximation and the numerical solution of PDEs are fundamental problems in scientific computing and applied mathematics \cite{quarteroni2006numerical,brenner2008mathematical,ciarlet2002finite,leveque2007finite}. Classical approaches, including polynomial and spline approximation, finite difference, finite element, finite volume, and discontinuous Galerkin methods, have established effective theoretical and computational frameworks for these problems \cite{leveque2002finite,arnold1982interior,reed1973triangular,riviere2008discontinuous}.
	
	% Function approximation and the numerical solution of PDEs are two fundamental problems in scientific computing and applied mathematics. Function approximation seeks to represent complex functions using accurate and computationally efficient representations \cite{quarteroni2006numerical}, while numerical methods for PDEs aim to obtain reliable approximations to solutions of mathematical models arising in science and engineering \cite{brenner2008mathematical, ciarlet2002finite,leveque2007finite}. 
	% %These two problems are closely related, since the numerical approximation of a PDE solution can often be viewed as a function approximation problem subject to the governing differential equation, boundary conditions, and other physical constraints. 
	% A variety of classical approaches, including polynomial and spline approximation, have been developed to construct accurate and efficient representations of functions. Similarly, the numerical solution of PDEs has led to a rich collection of discretization and approximation techniques, including finite difference \cite{leveque2007finite}, finite element \cite{brenner2008mathematical, ciarlet2002finite}, finite volume \cite{leveque2002finite}, and discontinuous Galerkin methods \cite{arnold1982interior, reed1973triangular, riviere2008discontinuous}. These methods have achieved remarkable success in a broad range of applications and have established well-developed theoretical and computational frameworks for function and PDE approximation.
	
	Recently, neural networks have emerged as powerful tools for function approximation and PDE computation due to their expressive power and flexibility. Beyond approximating general classes of functions, neural networks can also be used to represent numerical solution functions arising from conventional discretization methods, such as finite element approximations \cite{hao2026neural}. A variety of architectures, including multilayer perceptrons (MLPs), radial basis function (RBF) networks, sinusoidal representation networks (SIRENs), and Kolmogorov--Arnold networks (KANs), have been investigated for function approximation \cite{cybenko1989approximation,powell1987radial,sitzmann2020implicit,liu2025kan}. Neural networks have also been incorporated into PDE solvers, leading to methods such as physics-informed neural networks (PINNs), deep Ritz methods, Galerkin neural networks, neural Galerkin methods, and deep collocation methods \cite{raissi2019physics,yu2018deep,mark2021,mark2025,bruna2024neural,weng2026}. Most of these approaches rely on iterative optimization of network parameters, and their performance therefore depends critically on the training process.

	To reduce the cost of iterative training, non-iterative neural-network-based methods such as the extreme learning machine (ELM) have been developed \cite{huang2006universal,huang2006extreme}. In ELMs, hidden-layer parameters are randomly generated and fixed, while the output weights are determined by solving a least-squares problem. Originally developed for single-hidden-layer networks, ELMs have subsequently been extended to deeper architectures \cite{qu2016two,vong2018empirical,vasquez2023review}. Physics-informed extreme learning machines (PIELMs) further incorporate governing equations and boundary or initial conditions into the ELM framework, enabling efficient numerical solution of PDEs \cite{dwivedi2020physics,fabiani2021numerical,fabiani2025randonets}. However, ELM-type methods remain sensitive to the initialization of hidden-layer parameters \cite{dong2022computing,calabro2021extreme}.
	
	% To alleviate the computational cost associated with iterative training, non-iterative neural-network-based approximation strategies have also been developed. A notable example is the extreme learning machine (ELM) \cite{huang2006universal, huang2006extreme}, in which the hidden-layer parameters are randomly generated and fixed, while only the output-layer parameters are determined through a non-iterative procedure, typically by solving a least-squares problem. The ELM was originally proposed for single-hidden-layer feedforward neural networks \cite{huang2006extreme} and was subsequently extended to two-hidden-layer \cite{qu2016two} and multilayer architectures \cite{vong2018empirical}. A comprehensive review of multilayer ELMs and their variants can be found in \cite{vasquez2023review}, along with the references therein.
	% Building on ELM, physics-informed extreme learning machines (PIELMs) \cite{dwivedi2020physics} adopt the physics-informed learning framework of PINNs and incorporate the governing differential equations together with the associated boundary or initial conditions into the ELM framework. This combination enables the numerical solution of PDEs while retaining the non-iterative training strategy of ELMs, avoiding iterative optimization of the hidden-layer parameters. Subsequent developments have extended PIELMs to nonlinear PDEs, where perturbation techniques and nonlinear least-squares solvers are employed to address the resulting nonlinear algebraic systems \cite{fabiani2021numerical, fabiani2025randonets}.
	
	To improve the robustness of ELM-based PDE solvers, the Rank-Inspired Neural Network (RINN) framework was proposed in \cite{PengHuangYi2026}. Inspired by the orthogonality principle in spectral methods, RINN employs a covariance-based orthogonalization strategy to construct more stable neural basis functions. Its computation contains two stages: the neural basis is first trained and orthogonalized, after which the output weights are determined through linear or nonlinear least-squares solvers. Although RINN improves the stability of the approximation, its fully connected architecture can incur substantial computational costs, particularly for large networks and high-dimensional problems. Moreover, its performance is sensitive to a penalty parameter in the loss function, and the orthogonalization strategy does not explicitly promote a progressive increase in the effective dimensionality of the approximation space.
	
	% Although ELM-type methods offer significant computational efficiency advantages, their performance remains highly sensitive to the initialization of the hidden-layer parameters \cite{dong2022computing, calabro2021extreme}. 
	% %To address this issue, Dong et al. systematically investigated the effects of initialization and proposed a hyperparameter optimization framework based on residual analysis and differential evolution \cite{dong2022computing}. Calabr{\`o} et al. established theoretically justified parameter ranges for the hidden layers of PIELMs for one-dimensional elliptic problems with sharp boundary layers \cite{calabro2021extreme}. 
	% Inspired by the orthogonality principle underlying spectral methods, Peng et al. proposed the Rank-Inspired Neural Network (RINN) framework for solving PDEs \cite{PengHuangYi2026}, which introduces a covariance-based orthogonalization stage to mitigate the performance instability of PIELMs caused by random initialization. The RINN framework decouples the computational procedure into two stages: first, the neural basis functions are orthogonalized to improve the numerical stability of the approximation space; subsequently, the target function is approximated or the PDE is solved by determining the output weights through a least-squares problem, using linear or nonlinear least-squares solvers for linear or nonlinear problems, respectively. This strategy substantially improves the robustness of the resulting approximation with respect to the initialization parameters.
	
	In this work, we propose SCRINN for function approximation and the numerical solution of linear and nonlinear PDEs. Inspired by the localized basis functions used in the finite element method, SCRINN introduces sparse connectivity among hidden layers, such that each neural basis function depends only on a subset of neurons. This architecture reduces the number of network parameters and redundancy among neural basis functions, thereby improving computational efficiency. In addition, we propose a new robust loss function that replaces the logarithmic-type normalization used in RINN with a mean-square-type normalization. The resulting loss function substantially reduces the sensitivity to the normalization parameter and can also improve the performance of the original RINN; we refer to this variant as RINN+.

	Similar to RINN, SCRINN employs a two-stage framework. In the first stage, the neural basis functions are trained using the proposed loss function to promote orthogonality and normalization and construct an effective approximation space. In the second stage, the output weights are determined through a least-squares procedure. Extensive numerical experiments demonstrate that SCRINN achieves a favorable balance between accuracy, computational efficiency, and robustness across function approximation problems and steady-state, time-dependent, linear, and nonlinear PDEs. In particular, the sparse architecture significantly reduces computational costs while often achieving improved accuracy with a lower effective rank.

	The remainder of this paper is organized as follows. \Cref{subsecSCRINN} reviews the RINN algorithm and discusses its limitations. In \Cref{subsecRINN}, we introduce the proposed SCRINN method, including its sparse architecture, two-stage computational procedure, and robust loss function. Numerical experiments for function approximation and linear and nonlinear PDEs are presented in \Cref{sec:examples} and Supplementary Materials. Finally, \Cref{conclusion} summarizes the main conclusions and discusses possible future directions.
	
	% The remainder of this paper is organized as follows. \Cref{subsecSCRINN} reviews the RINN algorithm and discusses its underlying construction, including the generation of neural basis functions and the determination of the output weights. The limitations of the RINN framework are also examined. In \Cref{subsecRINN}, we introduce the SCRINN method and present its network architecture, sparse connectivity structure, and two-stage computational procedure. In particular, the proposed loss function for training the neural basis functions is described in detail. In \Cref{sec:examples}, we demonstrate the effectiveness of the proposed method through a series of numerical experiments, including function approximation, linear PDEs, and nonlinear PDEs. Finally, conclusions are drawn in \Cref{conclusion} to summarize the main findings and outline potential directions for future research.

	\section{Governing problems and RINN}\label{subsecSCRINN}
	In this work, we develop highly accurate deep learning models for two fundamental problems in numerical analysis: function approximation and the numerical solution of PDEs.
	
	\subsection{Governing problems}
	
	We consider the following two problems.
	
	\begin{problem}[Function approximation]
		\label{prob:function}
		Given a target function $u(\mathbf{x})$ defined on a domain
		$\mathcal{D}\subset\mathbb{R}^m$, the objective is to construct an accurate
		approximation $\hat{u}(\mathbf{x})$ such that
		\[
		\hat{u}(\mathbf{x})\approx u(\mathbf{x}),
		\qquad \mathbf{x}\in\mathcal{D}.
		\]
	\end{problem}
	
	% \paragraph{Problem 1: Function approximation}
	% Given a target function $u(\mathbf{x})$ defined on a domain
	% $\mathcal{D}$, the goal is to construct an accurate approximation
	% $\hat{u}(\mathbf{x})\approx u(\mathbf{x})$.
	
	\begin{problem}[Numerical solution of PDEs]
		\label{prob:pde}
		We consider two
		prototypical classes: time-dependent PDEs and elliptic PDEs. Both can be
		formulated in the unified form
		\begin{equation}\label{model}
			\begin{aligned}
				\mathscr{A}u(\mathbf{x}) &= f(\mathbf{x}), && \text{in } \mathcal{D},\\
				\mathscr{B}u(\mathbf{x}) &= g(\mathbf{x}), && \text{on } \Gamma,
			\end{aligned}
		\end{equation}
		where $\mathcal{D}$ denotes the problem domain, which may be either a
		spatial domain or a space--time domain, and $\Gamma$ denotes the boundary
		on which the boundary conditions are prescribed. The functions
		$f(\mathbf{x})$ and $g(\mathbf{x})$ are given source and boundary data,
		respectively. The operator $\mathscr{A}$ represents the differential
		operator associated with the governing PDE, while $\mathscr{B}$ denotes
		the boundary operator corresponding to the prescribed boundary
		conditions.

		%\paragraph{Problem 2a: Time-dependent PDEs}
		\noindent\textbf{Time-dependent PDEs.}
		For time-dependent problems, we consider the space--time domain
		\(
		\mathcal{D}=\Omega\times I,
		\)
		where $\Omega\subset\mathbb{R}^d$ is the spatial domain and
		$I=[0,T]$ with $T>0$ is the time domain. We write
		$\mathbf{x}=(x,t)$, where $x\in\Omega$ and $t\in I$.
		The boundary conditions are prescribed on the spatial boundary
		\(
		\Gamma=\partial\Omega\times I,
		\)
		while initial conditions are imposed on the initial-time surface
		\(
		\Gamma_0=\Omega\times\{0\}.
		\)
		To accommodate PDEs with different orders of time derivatives, we
		introduce the general initial-condition operator
		\begin{equation}\label{initial}
			\mathscr{I}u(\mathbf{x})=\mathbf{h}(\mathbf{x})
			\qquad \text{on } \Gamma_0,
		\end{equation}
		where $\mathscr{I}$ may represent one or multiple initial conditions, and $\mathbf{h}$ denotes the corresponding prescribed initial data.
		If the highest order of the time derivative is $r$, the required initial
		conditions \eqref{initial} can generally be written as
		\[
		\partial_t^k u(x,0)=h_k(x),
		\qquad k=0,\ldots,r-1,\quad x\in\Omega.
		\]
		For example, a first-order-in-time equation requires
		\[
		u(x,0)=h_0(x),
		\]
		whereas a second-order-in-time equation, such as the wave equation,
		requires
		\[
		u(x,0)=h_0(x),\qquad
		u_t(x,0)=h_1(x).
		\]
		
		\noindent\textbf{Elliptic PDEs.}
		For elliptic problems, the formulation is purely spatial, with
		$\mathcal{D} = \Omega$, $\mathbf{x} = x$, and $\Gamma = \partial\Omega$.
		Thus, \eqref{model} represents a standard boundary value problem, and
		no initial-condition operator is required.
	\end{problem}
	
	For notational convenience, we use the unified coordinate vector
	\[
	\mathbf{x}=(x_1,\ldots,x_m)\in\mathbb{R}^m,
	\]
	where the dimension
	\begin{equation}\label{dim}
		m=
		\begin{cases}
			d+1, & \text{for time-dependent PDEs},\\
			d, & \text{for elliptic PDEs}.
		\end{cases}
	\end{equation}

	\subsection{Overview of RINN}
	We first briefly review RINN, originally proposed in \cite{PengHuangYi2026} for solving PDEs, and then discuss its practical limitations.
	
	\subsubsection{Neural networks and neural basis functions}
	
	To approximate $u(\mathbf{x})$ in \eqref{model}, we introduce a fully connected neural network with $L$ hidden layers,
	\[
	u_{\theta}:\mathbb{R}^m\rightarrow\mathbb{R},
	\]
	defined by
	\begin{equation}\label{fullycon}
		u_{\theta}(\mathbf{x})
		=
		\left(
		\mathcal{A}_{\theta_{L+1}}
		\circ
		(\sigma\circ\mathcal{A}_{\theta_L})
		\circ\cdots\circ
		(\sigma\circ\mathcal{A}_{\theta_1})
		\right)(\mathbf{x}),
	\end{equation}
	where $\sigma:\mathbb{R}\rightarrow\mathbb{R}$ is the activation function,
	$\theta=(\theta_1,\ldots,\theta_{L+1})$ denotes the network parameters with
	$\theta_\ell=(W^{(\ell)},b^{(\ell)})$, and
	\[
	\mathcal{A}_{\theta_\ell}(z)
	=
	W^{(\ell)}z+b^{(\ell)},
	\qquad
	W^{(\ell)}\in\mathbb{R}^{N_\ell\times N_{\ell-1}},
	\quad
	b^{(\ell)}\in\mathbb{R}^{N_\ell}.
	\]
	
	Denote by $\phi_j(\hat{\theta},\mathbf{x})$ the output of the $j$-th neuron
	in the final hidden layer,
	\begin{equation}\label{nbf2}
		\phi_j(\hat{\theta},\mathbf{x})
		=
		\left[
		(\sigma\circ\mathcal{A}_{\theta_L})
		\circ\cdots\circ
		(\sigma\circ\mathcal{A}_{\theta_1})
		\right]_j(\mathbf{x}),
		\qquad
		1\leq j\leq N_L,
	\end{equation}
	where $\hat{\theta}=(\theta_1,\ldots,\theta_L)$. Taking
	$b^{(L+1)}=0$, the network output reduces to
	\begin{equation}\label{utheta}
		u_{\theta}(\mathbf{x})
		=
		\sum_{j=1}^{N_L}
		W^{(L+1)}_j\phi_j(\hat{\theta},\mathbf{x})
		=
		(W^{(L+1)})^\top
		\boldsymbol{\phi}(\hat{\theta},\mathbf{x}),
	\end{equation}
	where
	\[
	\boldsymbol{\phi}(\hat{\theta},\mathbf{x})
	=
	[\phi_1(\hat{\theta},\mathbf{x}),\ldots,
	\phi_{N_L}(\hat{\theta},\mathbf{x})]^\top
	\]
	is the vector of neural basis functions.
	
	To quantify the richness of these basis functions over the sampled domain,
	we define the neural basis matrix
	\begin{equation}\label{matrank}
		\mathbf{\Phi}
		=
		\begin{pmatrix}
			\boldsymbol{\phi}(\hat{\theta},\mathbf{x}_1)\\
			\boldsymbol{\phi}(\hat{\theta},\mathbf{x}_2)\\
			\vdots\\
			\boldsymbol{\phi}(\hat{\theta},\mathbf{x}_K)
		\end{pmatrix} = \begin{pmatrix}
			\phi _1( \hat{\theta} ,\mathbf{x}_1 )&   \phi _2( \hat{\theta} ,\mathbf{x}_1 )&   \cdots&   \phi _{N_L}( \hat{\theta}, \mathbf{x}_1 )\\
			\phi _1( \hat{\theta} ,\mathbf{x}_2 )&   \phi _2(  \hat{\theta} ,\mathbf{x}_2 )&   \cdots&   \phi _{N_L}( \hat{\theta} ,\mathbf{x}_2 )\\
			\vdots&   \vdots&   \ddots&   \vdots\\
			\phi _1(  \hat{\theta} , \mathbf{x}_K )&   \phi _2(  \hat{\theta} ,\mathbf{x}_K )&   \cdots&   \phi _{N_L}(  \hat{\theta} , \mathbf{x}_K )\\
		\end{pmatrix}
		\in\mathbb{R}^{K\times N_L},
	\end{equation}
	where $\{\mathbf{x}_i\}_{i=1}^K$ is the set of sample points with
	$K\geq N_L$. The rank of $\mathbf{\Phi}$ measures the linear independence
	of the neural basis functions on the sample set. A higher rank therefore
	indicates a richer neural basis, whereas rank deficiency indicates
	redundancy and may limit the approximation capability of the resulting
	neural representation.

	\begin{figure}[htbp]
		\centering
		\begin{subfigure}[b]{0.4\linewidth}
			\centering
			\includegraphics[width=1.0\linewidth]{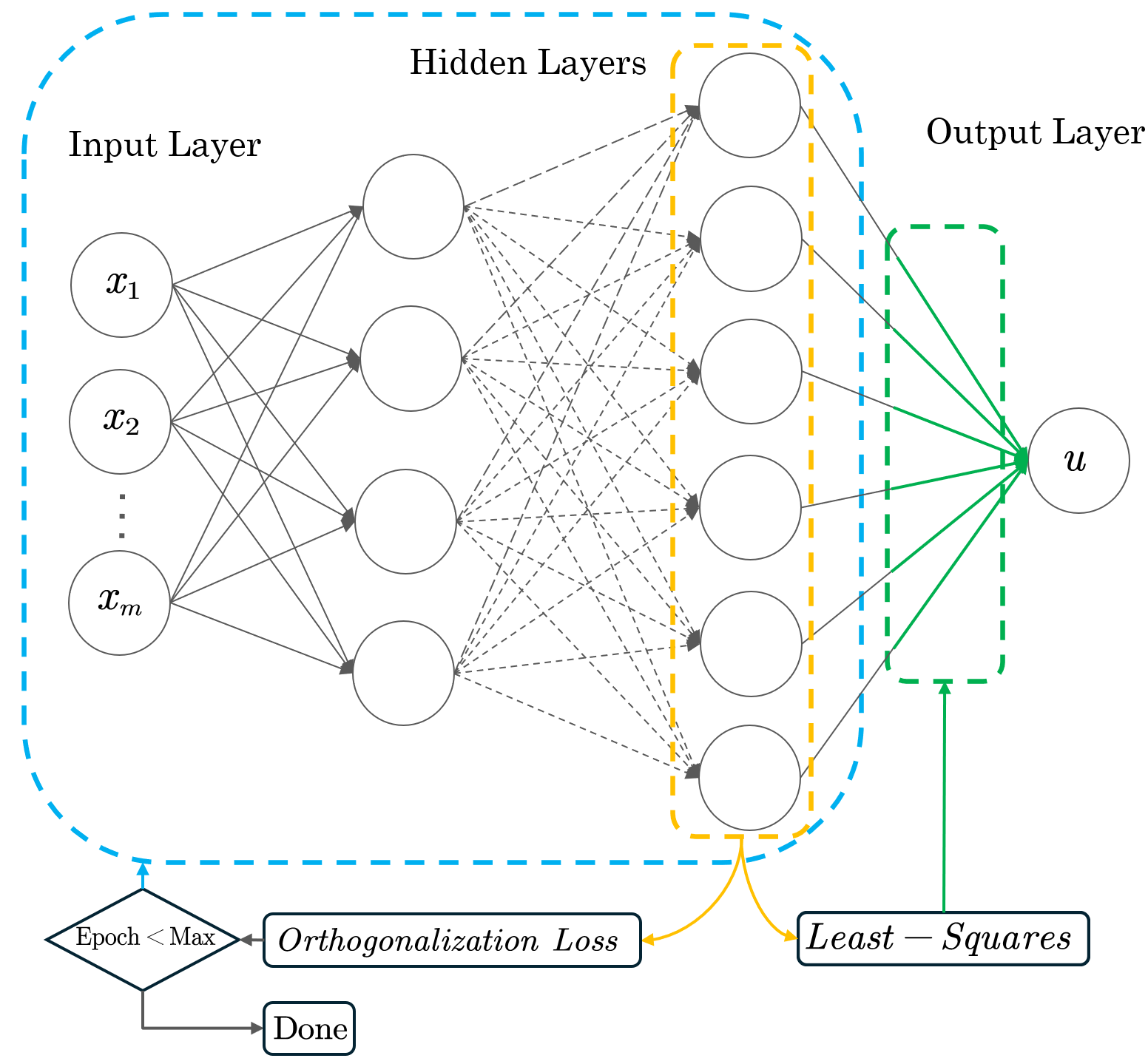}
			\caption{}
			\label{fig:rinn}
		\end{subfigure}
		\hspace{0.5cm}
		\begin{subfigure}[b]{0.4\linewidth}
			\centering
			\includegraphics[width=1.0\linewidth]{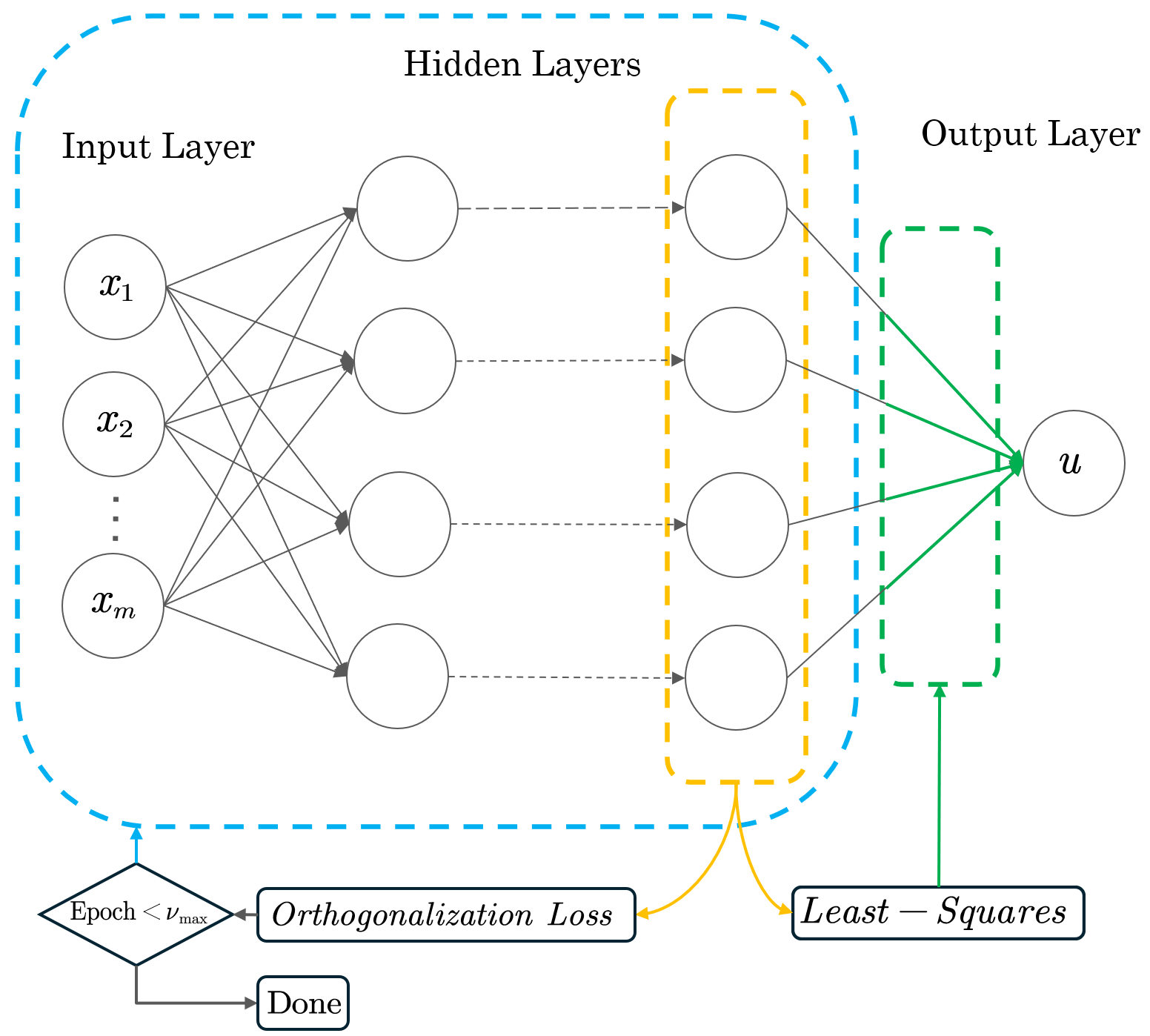}
			\caption{}
			\label{fig:scrinn}
		\end{subfigure}
		%  \caption{Neural network architectures. The blue dashed lines represent the trainable parameters $\hat{\theta}$ associated with the input and hidden layers, the yellow dashed lines denote the neural basis functions $\boldsymbol{\phi}(\hat{\theta},\mathbf{x})$ in the final hidden layer, and the green dashed lines represent the output-layer weights $W^{(L+1)}$ computed via a least-squares procedure. (a) RINN; (b) SCRINN.}
		\caption{Neural network architectures: (a) RINN and (b) SCRINN. Blue, yellow, and green dashed lines denote the trainable parameters $\hat{\theta}$, neural basis functions $\boldsymbol{\phi}(\hat{\theta},\mathbf{x})$, and least-squares output weights $W^{(L+1)}$, respectively.}\label{fig:nna}
	\end{figure}
	
	\subsubsection{RINN framework}
	RINN was developed in \cite{PengHuangYi2026} to address the sensitivity of PIELMs to weight initialization. The architecture of RINN is illustrated in \Cref{fig:rinn}. 
	%The blue dashed lines represent the trainable parameters $\hat{\theta}$ associated with the input and hidden layers, the yellow dashed lines denote the neural basis functions $\boldsymbol{\phi}(\hat{\theta},\mathbf{x})$ in the final hidden layer, and the green dashed lines represent the output-layer weights $W^{(L+1)}$. %, which are determined through least-squares fitting.
	% \begin{figure}[htbp]
		%   \centering
		%   \includegraphics[width=0.4\linewidth]{rinn}\hspace{1.5cm}
		%   \includegraphics[width=0.4\linewidth]{srinn}
		%   \caption{Neural network architectures. Left: RINN, Right: SCRINN}\label{fig:nna}
		% \end{figure}
	RINN is formulated as a two-stage framework. In the first stage, the hidden-layer parameters $\hat{\theta}$ are optimized so that the resulting neural basis functions satisfy an approximate orthogonality condition. Ideally, the neural basis functions are orthogonal with respect to the $L^2$-inner product,
	\begin{equation}\label{L2cond1}
		(\phi_i,\phi_j)_{L^2(\mathcal{D})}
		=
		\int_{\mathcal{D}}
		\phi_i(\hat{\theta},\mathbf{x})\,\phi_j(\hat{\theta},\mathbf{x})
		\, d\mathbf{x}
		= 0,
		\qquad i\neq j.
	\end{equation}
	while each basis function is normalized in the $L^2$ norm, i.e.,
	\begin{equation}\label{L2normal}
		\| \phi _i\|_{L^{2}\left( \mathcal{D} \right)}^2 =
		\int_{\mathcal{D}}
		\phi_i(\hat{\theta},\mathbf{x})^2
		\, d\mathbf{x}
		= 1, \quad \forall i.
	\end{equation}
	In practice, the orthogonality \eqref{L2cond1} is approximated using a set of collocation points in domain $\mathcal{D}$. This yields the discrete orthogonality constraint on 
	$\{\mathbf{x}_k\}_{k=1}^{K}\subset\mathcal{D}$ by %the discrete $L^2$-inner product,
	\begin{equation}\label{normalized}
		\sum_{k=1}^K{ \phi _i( \hat{\theta} ,\mathbf{x}_k ) \cdot \,\phi _j( \hat{\theta} ,\mathbf{x}_k )} =0.
	\end{equation}
	
	To train the neural basis functions so that they approximately satisfy the
	orthogonality and normalization conditions \eqref{L2cond1} and
	\eqref{L2normal}, respectively, in a discrete sense, the sample covariance matrix of $\mathbf{\Phi}$ was introduced:
	\begin{equation}\label{matrixC1}
		\mathbf{C}=\frac{1}{K-1}\mathbf{\Phi}^{\top}\mathbf{\Phi},
	\end{equation}
	where the entry
	\begin{equation}
		\mathbf{C}_{ij} = \frac{1}{K-1} \sum_{k=1}^K{ \phi _i( \hat{\theta} ,\mathbf{x}_k ) \cdot \,\phi _j( \hat{\theta} ,\mathbf{x}_k )}
	\end{equation}
	provides a scaled discrete approximation of the $L^2$ inner product
	between the neural basis functions $\phi_i$ and $\phi_j$. 
	Accordingly, the orthogonality condition \eqref{L2cond1} and normalization
	condition \eqref{L2normal} can be approximately enforced in the discrete
	sense by driving the off-diagonal and diagonal entries of $\mathbf{C}$
	toward zero and one, respectively. This is achieved by minimizing the total loss
	\begin{equation}\label{lossf}
		\mathcal{L}_{\mathrm{RINN}}(\hat\theta) =  \mathcal{L}_\mathrm{orthogonal}(\hat\theta) + \varepsilon \mathcal{L}_\mathrm{normal}(\hat\theta),
	\end{equation}
	where $\varepsilon>0$ is a tunable parameter that controls the relative weight of the normalization term. The value $\varepsilon = 0.01$ is chosen empirically \cite{PengHuangYi2026}.
	The orthogonality loss is defined as
	\begin{equation}\label{lossortho}
		\mathcal{L}_\mathrm{orthogonal} (\hat\theta) =\sqrt{\sum\limits_{ 1\leq  i, j\leq N_L, i\ne j }|\mathbf{C}_{ij}|^2},
	\end{equation}
	which corresponds to the Frobenius norm of the off-diagonal entries of $\mathbf{C}$.
	%Minimizing $\mathcal{L}_\mathrm{orthogonal}$ drives the off-diagonal entries of $\mathbf{C}$ toward zero, thereby enforcing the discrete orthogonality condition \eqref{L2cond1}.
	The normalization loss is defined by
	\begin{equation}\label{lossnormal}
		\mathcal{L}_\mathrm{normal}(\hat\theta) =\sum\limits_{i=1}^{N_L}\left | \log_{10}(\mathbf{C}_{ii}^2)\right |=\sum\limits_{i=1}^{N_L}\left | 2\log_{10}(\mathbf{C}_{ii})\right |,
	\end{equation}
	which penalizes deviations of the diagonal entries of $\mathbf{C}$ from one. 
	
	In the second stage, the hidden-layer parameters $\hat{\theta}$ are kept fixed, and the output-layer parameters $W^{(L+1)}$ are computed using a least-squares method to approximate the solution of the PDEs at selected collocation points.
	
	%To improve the robustness of RINN, an early stopping strategy can be incorporated during the first-stage training, resulting in the variant RINN-es. Numerical tests in \cite{PengHuangYi2026} demonstrate that, compared with PIELM, RINN-es significantly reduces sensitivity to parameter initialization while maintaining high solution accuracy.
	
	To enhance the robustness of RINN, \cite{PengHuangYi2026} introduced an
	early stopping strategy for the first-stage training, resulting in the
	RINN-es variant. Numerical experiments demonstrate that RINN-es
	substantially reduces the sensitivity to parameter initialization compared
	with PIELM, while maintaining high solution accuracy.

	\subsection{Limitations of RINN}
	Despite the advantages, numerical investigations reveal several limitations of the RINN framework. First, RINN employs a fully connected neural network architecture with a large number of parameters, which leads to substantial computational and optimization costs in the initial training stage. In particular, for high-dimensional PDEs, the training procedure becomes increasingly infeasible as the number of collocation points grows.
	Second, the tunable parameter $\varepsilon$ in \eqref{lossf} should be chosen small, which is far less than $1$, and its optimal value may depend on the problem, while larger values typically result in larger errors.
	Third, while the sample covariance matrix is designed to enforce the orthogonality condition \eqref{L2cond1} and normalization condition \eqref{L2normal}, it does not guarantee progressively stronger linear independence among the learned neural basis functions.
	As a consequence, the expressive power of the induced approximation space is not guaranteed to improve monotonically. From the perspective of the neural basis function matrix $\mathbf{\Phi}$ in \eqref{matrank}, this implies that its rank is not guaranteed to increase throughout training and may even decrease, thereby degrading the approximation capability.
	
	% The first issue is readily observed in practice due to the rapid growth of computational cost and the associated optimization instability. 
	% For the second issue, $\varepsilon=0.01$ was recommended for many examples in \cite{PengHuangYi2026}, while larger values may introduce greater approximation errors. The third issue can be viewed from the rank evolution of RINN-es, see \Cref{fig:posong3.3.2} in \Cref{ex:2DPoisson}.
	
	The first issue is readily observed in practice, as increasing the network size leads to rapidly growing computational costs and greater optimization instability. For the second issue, $\varepsilon=0.01$ was recommended for many examples in \cite{PengHuangYi2026}, whereas larger values may result in greater approximation errors. The third issue is reflected in the rank evolution of RINN-es; see \Cref{fig:posong3.3.2} in \Cref{ex:2DPoisson}.

	\section{SCRINN}\label{subsecRINN}
	To overcome the limitations of RINN discussed above, we propose SCRINN, a near-orthogonalization training strategy enabled by a sparse connectivity architecture. The key objective of SCRINN is 
	to reduce the parameter coupling among neural basis functions, thereby mitigating rank saturation and facilitating sustained growth of the rank of the neural basis function matrix $\mathbf{\Phi}$ in \eqref{matrank}.
	
	\subsection{Sparsely connected neural networks}

	Unlike the fully connected architecture used in RINN (see \Cref{fig:rinn}), 
	SCRINN employs a sparsely connected architecture, as illustrated in 
	\Cref{fig:scrinn}. The input and output layers remain fully connected, 
	while the intermediate hidden layers adopt a channel-wise sparse 
	connectivity pattern. Let $N_1=\cdots=N_L=N$ denote the number of 
	neurons in the hidden layers. The resulting network is defined as
	\begin{equation} \label{notfullycon}
		u_{\theta}(\mathbf{x})
		=
		\left(
		\mathcal{S}_{\theta_{L+1}}
		\circ
		(\sigma\circ\mathcal{S}_{\theta_L})
		\circ\cdots\circ
		(\sigma\circ\mathcal{S}_{\theta_2})
		\circ
		(\sigma\circ\mathcal{S}_{\theta_1})
		\right)(\mathbf{x}),
	\end{equation}
	where $\sigma:\mathbb{R}\rightarrow\mathbb{R}$ is the activation function
	and $\theta=(\theta_1,\ldots,\theta_{L+1})$ denotes the network parameters,
	with $\theta_\ell=(\hat W^{(\ell)},b^{(\ell)})$. For $\ell=1 $, and $ \ell= L+1$, $\mathcal{S}_{\theta_\ell}$ is the standard affine
	transformation
	\[
	\mathcal{S}_{\theta_\ell}(z)
	=
	\hat W^{(\ell)}z+b^{(\ell)}, \qquad
	\hat W^{(\ell)}\in\mathbb{R}^{N_\ell\times N_{\ell-1}},
	\quad
	b^{(\ell)}\in\mathbb{R}^{N_\ell}.
	\]
	whereas for $2\leq\ell\leq L$,
	\[
	\mathcal{S}_{\theta_\ell}(z)
	=
	z\odot\hat W^{(\ell)}+b^{(\ell)},
	\qquad
	\hat W^{(\ell)},b^{(\ell)}\in\mathbb{R}^{N}.
	\]
	Here, $\odot$ denotes the Hadamard product, corresponding to diagonal
	connectivity in the intermediate hidden layers.
	
	Denote by $\phi_j(\hat{\theta},\mathbf{x})$ the output of the $j$-th neuron
	in the final hidden layer,
	\begin{equation}\label{nbf}
		\phi_j(\hat{\theta},\mathbf{x})
		=
		\left[
		(\sigma\circ\mathcal{S}_{\theta_L})
		\circ\cdots\circ
		(\sigma\circ\mathcal{S}_{\theta_1})
		\right]_j(\mathbf{x}),
		\qquad
		1\leq j\leq N,
	\end{equation}
	where $\hat{\theta}=(\theta_1,\ldots,\theta_L)$. With $b^{(L+1)}=0$,
	the network output can then be written as
	\begin{equation}\label{solscrinn}
		u_{\theta}(\mathbf{x})
		=
		\sum_{j=1}^{N}
		\hat W^{(L+1)}_j\phi_j(\hat{\theta},\mathbf{x})
		=
		(\hat W^{(L+1)})^\top
		\boldsymbol{\phi}(\hat{\theta},\mathbf{x}),
	\end{equation}
	where
	\begin{equation}\label{basisvector}
		\boldsymbol{\phi}(\hat{\theta},\mathbf{x})
		=
		[\phi_1(\hat{\theta},\mathbf{x}),\ldots,\phi_N(\hat{\theta},\mathbf{x})]^\top
	\end{equation}
	denotes the vector of neural basis functions. The corresponding matrix
	representation can be given in the form of \eqref{matrank}.
	
	% \begin{remark}
		%     The sparsely connected neural network based on the architecture in \Cref{fig:scrinn} can also be interpreted as a special case of a fully connected neural network (or residual neural network), in which all weight matrices, except those associated with the input and output layers, are constrained to be diagonal.
		% It is of interest to investigate more general sparsely connected architectures, for instance, where different hidden layers contain different numbers of neurons, and each neuron is connected to a fixed subset of neurons in adjacent layers rather than exactly one or all neurons. Such structured sparsity may provide a balance between model expressivity and computational efficiency. This direction will be pursued in future work.
		% \end{remark}
	
	\begin{remark}
		The sparsely connected architecture in \Cref{fig:scrinn} can be viewed as a special case of a fully connected neural network where the weight matrices of the intermediate hidden layers are constrained to be diagonal. More general sparse architectures, allowing different numbers of neurons across hidden layers and connections to a fixed subset of neighboring neurons, may provide a favorable balance between expressivity and computational efficiency and will be investigated in future work.
	\end{remark}

	%\subsection{SCRINN framework}
	% The idea of SCRINN is inspired by the finite element method, where compactly supported basis functions are employed to construct localized and numerically stable approximations. Analogously, SCRINN introduces sparse connectivity among neurons in hidden layers, such that each neural basis function depends only on a small subset of hidden-layer neurons. This ``localized" dependency structure reduces redundancy among neural basis functions and promotes an increase in the rank of the neural basis function matrix during training, thereby effectively improving solution accuracy. The proposed SCRINN follows a similar two-stage framework as RINN. 
	
	\subsection{SCRINN framework}
	
	Inspired by the finite element method, where localized basis functions provide effective and stable approximations, SCRINN introduces sparse connectivity among hidden-layer neurons to reduce parameter coupling and redundancy among neural basis functions. This structured sparsity is designed to promote the growth of the neural basis matrix rank and thereby enhance the approximation capability. The proposed SCRINN also follows a two-stage framework.
	
	\subsubsection{Stage I: Orthogonality and normalization of neural basis functions}
	In the first stage, the neural basis functions are trained to satisfy the $L^2$-orthogonality condition \eqref{L2cond1} and the normalization condition \eqref{L2normal} by minimizing a loss function analogous to \eqref{lossf}. However, SCRINN employs a different approximation of the $L^2$ inner product and a different formulation of the loss function. The orthogonalization process serves to reduce redundancy among the neural basis functions and promote their linear independence, thereby improving the rank of the resulting neural basis function matrix, as quantified in \eqref{matrank}.
	
	To approximate the $L^2$ inner products in \eqref{L2cond1} and \eqref{L2normal}, we introduce two numerical integrations.
	The first approach is based on Monte Carlo integration. Let
	$\{\mathbf{x}_k\}_{k=1}^{K}\subset\mathcal{D}$ be a set of collocation points with
	$K\ge N_L$. The discrete $L^2$ inner product is approximated by
	\begin{equation}\label{discreteH1}
		\left(\phi_i,\phi_j\right)_{L_h^2(\mathcal D)}
		=
		\frac{|\mathcal D|}{K}
		\sum_{k=1}^{K}
		\phi_i(\hat\theta,\mathbf{x}_k)
		\phi_j(\hat\theta,\mathbf{x}_k),
	\end{equation}
	where $|\mathcal D|$ denotes the $m$-dimensional measure (length, area, volume, or hypervolume) of the domain $\mathcal D$.
	The second approach employs Gaussian quadrature. 
	Let $\{(\mathbf{x}_k,\omega_k)\}_{k=1}^{K}$ denote the quadrature points and weights induced by a Gaussian quadrature rule applied on each element of a partition of domain $\mathcal{D}$. The discrete $L^2$ inner product is approximated by
	\begin{equation}\label{discreteGaussian}
		\left(\phi_i,\phi_j\right)_{L_h^2(\mathcal D)}
		=
		\sum_{k=1}^{K}
		\omega_k\,
		\phi_i(\hat\theta,\mathbf{x}_k)
		\phi_j(\hat\theta,\mathbf{x}_k).
	\end{equation}
	Unless otherwise specified, either \eqref{discreteH1} or \eqref{discreteGaussian} is used to define the discrete $L^2$ inner product, with the associated discrete norm given by
	\[
	\|\phi_i\|_{L_h^2(\mathcal D)}
	=
	\sqrt{
		\left(\phi_i,\phi_i\right)_{L_h^2(\mathcal D)}
	}.
	\]
	
	Accordingly, the continuous $L^2$ orthogonality condition \eqref{L2cond1} is replaced by the discrete $L^2$ orthogonality condition
	\begin{equation}\label{2.3liner11}
		\left( \phi _i,\phi _j \right)_{L^2_h\left( \mathcal{D} \right)}
		= 0,
		\qquad i\neq j,
	\end{equation}
	while the normalization condition \eqref{L2normal} becomes
	\begin{equation}\label{discretenormal}
		\| \phi _i\|_{L^2_h\left( \mathcal{D} \right)}^2=1,
		\qquad \forall i.
	\end{equation}
	
	To enforce the discrete $L^2$ orthogonality and normalization conditions in \eqref{2.3liner11} and \eqref{discretenormal}, we introduce the Gram matrix
	\begin{equation}\label{2.3liner3}
		\mathbf{G}
		=
		\big[
		(\phi_i,\phi_j)_{L_h^2(\mathcal{D})}
		\big]_{i,j=1}^{N_L}.
	\end{equation}
	%Thus, the desired orthonormality condition is equivalent to $\mathbf{G}=\mathbf{I}$.
	% To enforce the discrete $L^2$ orthogonality condition \eqref{2.3liner11} and the discrete normalization condition \eqref{discretenormal}, we introduce the Gram matrix associated with the discrete $L_h^2$ inner product:
	% \begin{equation}\label{2.3liner3}
		% \mathbf{G}
		% =
		% \left[
		% \begin{matrix}
			% (\phi_1,\phi_1)_{L_h^2(\mathcal{D})}
			% &
			% (\phi_1,\phi_2)_{L_h^2(\mathcal{D})}
			% &
			% \cdots
			% &
			% (\phi_1,\phi_{N_L})_{L_h^2(\mathcal{D})}
			% \\
			% (\phi_2,\phi_1)_{L_h^2(\mathcal{D})}
			% &
			% (\phi_2,\phi_2)_{L_h^2(\mathcal{D})}
			% &
			% \cdots
			% &
			% (\phi_2,\phi_{N_L})_{L_h^2(\mathcal{D})}
			% \\
			% \vdots
			% &
			% \vdots
			% &
			% \ddots
			% &
			% \vdots
			% \\
			% (\phi_{N_L},\phi_1)_{L_h^2(\mathcal{D})}
			% &
			% (\phi_{N_L},\phi_2)_{L_h^2(\mathcal{D})}
			% &
			% \cdots
			% &
			% (\phi_{N_L},\phi_{N_L})_{L_h^2(\mathcal{D})}
			% \end{matrix}
		% \right].
		% \end{equation}
	% Then, the orthogonality and normalization conditions can be enforced by driving the off-diagonal entries of the Gram matrix $\mathbf{G}$ toward zero and its diagonal entries toward one, respectively. Since $\mathbf{G}$ is symmetric, it suffices to consider only its strictly upper triangular part when measuring orthogonality. Consequently, instead of the loss function \eqref{lossf}, we define the SCRINN loss as
	Since $\mathbf{G}$ is symmetric, the orthogonality and normalization
	conditions can be enforced by driving its strictly upper triangular
	entries to zero and its diagonal entries to one, respectively. Thus, we
	define the SCRINN loss by
	\begin{equation}\label{lossscrinn}
		\mathcal{L}_{\mathrm{SCRINN}}(\hat\theta) =  \mathcal{L}_\mathrm{orthogonal}(\hat\theta) + \varepsilon\mathcal{L}_\mathrm{normal}(\hat\theta),
	\end{equation}
	where $\varepsilon>0$ controls the relative weight of the normalization
	term. The two components are given by
	\begin{equation}\label{lossortho+}
		\mathcal{L}_{\mathrm{orthogonal}}(\hat\theta)
		=
		\sum_{1\leq i<j\leq N_L}|\mathbf{G}_{ij}|^2,
	\end{equation}
	and
	\begin{equation}\label{lossnormal+}
		\mathcal{L}_{\mathrm{normal}}(\hat\theta)
		=
		\sum_{i=1}^{N_L}|\mathbf{G}_{ii}-1|^2.
	\end{equation}
	The former penalizes the off-diagonal entries of $\mathbf{G}$ and hence
	promotes discrete orthogonality, while the latter drives the diagonal
	entries toward unity and enforces discrete normalization. Unlike the
	normalization strategy adopted in RINN, the discrete inner product in
	\eqref{discreteH1} directly approximates the $L^2$ inner product, which
	motivates the normalization loss in \eqref{lossnormal+}.

	\subsubsection{Stage II: Least-squares problem for function fitting and PDE solving}
	Let $\hat{\theta}^{(\nu)}$ denote the network parameters at
	the $\nu$-th epoch obtained by minimizing \eqref{lossscrinn}.
	In Stage II, the parameters $\hat{\theta}^{(\nu)}$ are fixed, so that the trained neural basis functions $\boldsymbol{\phi}(\hat{\theta}^{(\nu)},\mathbf{x})$ in \eqref{basisvector} remain unchanged. Only the output-layer weights are determined in this stage. The resulting procedure follows the standard ELM framework for function approximation and the PIELM framework for PDE solving. We briefly describe the formulation below, with particular attention to the treatment of nonlinear PDEs.
	
	To distinguish the learning strategies in the two stages, we denote the output-layer weight vector $W^{(L+1)}$ by
	\[
	\boldsymbol{\beta}
	=
	[\beta_1,\ldots,\beta_{N_L}]^\top.
	\]
	Accordingly, the network output takes the form
	\begin{equation}\label{uthetaII}
		u_{\theta}(\mathbf{x})
		=
		\sum_{j=1}^{N_L}\beta_j
		\phi_j(\hat{\theta}^{(\nu)},\mathbf{x})
		=
		\boldsymbol{\beta}^\top
		\boldsymbol{\phi}(\hat{\theta}^{(\nu)},\mathbf{x}).
	\end{equation}
	
	\noindent
	\textbf{Function approximation.}
	We consider function approximation using \eqref{uthetaII}, in which the neural basis functions $\boldsymbol{\phi}(\hat{\theta}^{(\nu)}, \mathbf{x})$ are fixed. 
	Let
	$\mathbf{X}\in\mathbb{R}^{K\times m}$ denote a set of collocation points in $\overline{\mathcal{D}}$, and let
	$u(\mathbf{X})\in\mathbb{R}^K$ denote the corresponding function values. Evaluating the fixed neural basis functions at these points gives the linear system
	\begin{equation}\label{2.2liner1approximation}
		\mathbf{H}(\hat{\theta}^{(\nu)},\mathbf{X})\boldsymbol{\beta}
		=
		\mathbf{S},
	\end{equation}
	where
	\[
	\mathbf{H}(\hat{\theta}^{(\nu)},\mathbf{X})
	=
	\boldsymbol{\phi}(\hat{\theta}^{(\nu)},\mathbf{X})
	\in\mathbb{R}^{K\times N_L},
	\qquad
	\mathbf{S}=u(\mathbf{X})\in\mathbb{R}^K.
	\]
	Since $K$ is typically larger than $N_L$, the system \eqref{2.2liner1approximation} is generally overdetermined. We therefore determine the output weights by solving the least-squares problem
	\[
	\boldsymbol{\beta}^{(\nu)}
	=
	\arg\min_{\boldsymbol{\beta}}
	\left\|
	\mathbf{H}(\hat{\theta}^{(\nu)},\mathbf{X})\boldsymbol{\beta}
	-
	\mathbf{S}
	\right\|_2^2.
	\]
	A minimum-norm least-squares solution is given by
	\begin{equation}\label{lspsol}
		\boldsymbol{\beta}^{(\nu)}
		=
		\mathbf{H}(\hat{\theta}^{(\nu)},\mathbf{X})^\dagger
		\mathbf{S},
	\end{equation}
	where $\mathbf{H}^\dagger$ denotes the Moore--Penrose pseudoinverse of $\mathbf{H}$. The resulting approximation at any query point $\mathbf{x}$ is then
	\begin{equation}\label{uthetaII+}
		u_{\theta}^{(\nu)}  (\mathbf{x})
		=
		(\boldsymbol{\beta}^{(\nu)})^\top
		\boldsymbol{\phi}(\hat{\theta}^{(\nu)},\mathbf{x}).
	\end{equation}
	
	\noindent
	\textbf{Linear PDE solving.}
	For a linear PDE, we substitute the expansion \eqref{uthetaII} with the frozen neural basis functions into the governing equation \eqref{model}. All terms involving the initial condition are included only for time-dependent problems. The output weights are then determined by a least-squares formulation based on collocation points sampled from the relevant parts of the computational domain.
	
	Specifically, we sample residual points
	$\mathbf{x}^{\mathrm{res}}\in\mathbb{R}^{K^{\mathrm{res}}\times m}$ in $\mathcal{D}$, boundary points
	$\mathbf{x}^{\mathrm{bcs}}\in\mathbb{R}^{K^{\mathrm{bcs}}\times m}$ on $\Gamma$, and, for time-dependent problems, initial points
	$\mathbf{x}^{\mathrm{ics}}\in\mathbb{R}^{K^{\mathrm{ics}}\times m}$ on $\Gamma_0$. We collect these points into the collocation set
	\[
	\mathbf{X}
	=
	\left\{
	\mathbf{x}^{\mathrm{res}},
	\mathbf{x}^{\mathrm{bcs}},
	\mathbf{x}^{\mathrm{ics}}
	\right\}
	\in
	\mathbb{R}^{K\times m},
	\]
	where
	\[
	K
	=
	K^{\mathrm{res}}
	+
	K^{\mathrm{bcs}}
	+
	K^{\mathrm{ics}}.
	\]
	For stationary problems, the initial-point set and the corresponding terms are simply omitted.

	Enforcing the governing equation, boundary conditions, and, for
	time-dependent problems, initial conditions at the corresponding
	collocation points yields the following linear system:
	\begin{equation}\label{2.3out2}
		\left\{
		\begin{aligned}
			&\sum_{i=1}^{N_L}
			\left[\mathscr{A}\phi_i(\hat{\theta}^{(\nu)},\mathbf{x}^{\mathrm{res}})\right]\beta_i
			=
			f(\mathbf{x}^{\mathrm{res}}),
			&&\mathbf{x}^{\mathrm{res}}\in\mathcal{D},
			\\
			&\sum_{i=1}^{N_L}
			\left[\mathscr{B}\phi_i(\hat{\theta}^{(\nu)},\mathbf{x}^{\mathrm{bcs}})\right]\beta_i
			=
			g(\mathbf{x}^{\mathrm{bcs}}),
			&&\mathbf{x}^{\mathrm{bcs}}\in\Gamma,
			\\
			&\sum_{i=1}^{N_L}
			\left[\mathscr{I}\phi_i(\hat{\theta}^{(\nu)},\mathbf{x}^{\mathrm{ics}})\right]\beta_i
			=
			\mathbf{h}(\mathbf{x}^{\mathrm{ics}}),
			&&\mathbf{x}^{\mathrm{ics}}\in\Gamma_0.
		\end{aligned}
		\right.
	\end{equation}
	Here, the terms involving the differential, boundary, and initial
	operators are evaluated using forward propagation and automatic
	differentiation.
	
	Since the neural basis functions are fixed in Stage II, the above
	system is linear with respect to the output weights $\boldsymbol{\beta}$.
	It can therefore be written compactly as
	\begin{equation}\label{2.2liner14}
		\mathbf{H}(\hat{\theta}^{(\nu)},\mathbf{X})
		\boldsymbol{\beta}
		=
		\mathbf{S},
	\end{equation}
	where
	\[
	\mathbf{H}(\hat{\theta}^{(\nu)},\mathbf{X})
	=
	\begin{bmatrix}
		\mathscr{A}\boldsymbol{\phi}(\hat{\theta}^{(\nu)},\mathbf{x}^{\mathrm{res}})
		\\
		\mathscr{B}\boldsymbol{\phi}(\hat{\theta}^{(\nu)},\mathbf{x}^{\mathrm{bcs}})
		\\
		\mathscr{I}\boldsymbol{\phi}(\hat{\theta}^{(\nu)},\mathbf{x}^{\mathrm{ics}})
	\end{bmatrix}
	\in\mathbb{R}^{K\times N_L},
	\qquad
	\mathbf{S}
	=
	\begin{bmatrix}
		f(\mathbf{x}^{\mathrm{res}})
		\\
		g(\mathbf{x}^{\mathrm{bcs}})
		\\
		\mathbf{h}(\mathbf{x}^{\mathrm{ics}})
	\end{bmatrix}
	\in\mathbb{R}^{K}.
	\]
	The linear system \eqref{2.2liner14} is solved in the same manner as
	\eqref{2.2liner1approximation}, yielding the output weights
	$\boldsymbol{\beta}^{(\nu)}$. The resulting neural network solution is then
	given by \eqref{uthetaII+}.
	
	\noindent
	\textbf{Nonlinear PDE solving.}
	For nonlinear PDEs, the corresponding collocation system is generally
	nonlinear with respect to the output weights $\boldsymbol{\beta}$,
	because the differential operator may contain nonlinear functions or
	products of the solution. Therefore, the output weights are determined
	by solving a nonlinear least-squares problem.
	
	Substituting the neural representation \eqref{uthetaII} into the
	governing equation and the associated boundary and initial conditions
	gives the residual vector
	\begin{equation}\label{residualeqn}
		\mathbf{R}(\boldsymbol{\beta})
		=
		\begin{bmatrix}
			\mathscr{A}\bigl[u_{\theta}(\mathbf{x}^{\mathrm{res}})\bigr]
			-
			f(\mathbf{x}^{\mathrm{res}})
			\\
			\mathscr{B}\bigl[u_{\theta}(\mathbf{x}^{\mathrm{bcs}})\bigr]
			-
			g(\mathbf{x}^{\mathrm{bcs}})
			\\
			\mathscr{I}\bigl[u_{\theta}(\mathbf{x}^{\mathrm{ics}})\bigr]
			-
			\mathbf{h}(\mathbf{x}^{\mathrm{ics}})
		\end{bmatrix},
	\end{equation}
	where $u_{\theta}(\mathbf{x})$ has the form of \eqref{uthetaII}.
	% \[
	% u_{\theta}(\mathbf{x})
	% =
	% \boldsymbol{\beta}^{\top}
	% \boldsymbol{\phi}(\hat{\theta}^{(\nu)},\mathbf{x}).
	% \]
	For stationary problems, the third block corresponding to the initial
	condition is omitted. The output weights are then determined by solving
	the nonlinear least-squares problem
	\begin{equation}\label{nonlinearLS}
		\boldsymbol{\beta}^{(\nu)}
		=
		\arg\min_{\boldsymbol{\beta}}
		\left\|
		\mathbf{R}(\boldsymbol{\beta})
		\right\|_2^2.
	\end{equation}
	
	Unlike the linear PDE case, the residual vector
	$\mathbf{R}(\boldsymbol{\beta})$ is generally nonlinear in
	$\boldsymbol{\beta}$. Consequently, the direct pseudoinverse solution
	in \eqref{lspsol} is no longer applicable. Instead, we employ an
	iterative nonlinear least-squares solver. Specifically, the
	Trust Region Reflective (TRF) algorithm \cite{branch1999subspace} is
	used to update the output weights
	\[
	\boldsymbol{\beta}^{(k+1)}
	=
	\boldsymbol{\beta}^{(k)}
	+
	\Delta\boldsymbol{\beta}^{(k)},
	\]
	where $\Delta\boldsymbol{\beta}^{(k)}$ is obtained by solving the
	trust-region subproblem
	\begin{equation}\label{TRsubproblem}
		\Delta\boldsymbol{\beta}^{(k)}
		=
		\arg\min_{\|\Delta\boldsymbol{\beta}\|_2\leq \Delta_k}
		\left\|
		\mathbf{R}\bigl(\boldsymbol{\beta}^{(k)}\bigr)
		+
		\mathbf{J}\bigl(\boldsymbol{\beta}^{(k)}\bigr)
		\Delta\boldsymbol{\beta}
		\right\|_2^2,
	\end{equation}
	where $\Delta_k$ denotes the trust-region radius at iteration $k$.
	The radius is adaptively adjusted by the TRF algorithm based on the
	actual and predicted reductions of the objective function. Here,
	\begin{equation}\label{jacobian}
		\mathbf{J}(\boldsymbol{\beta})
		=
		\frac{\partial\mathbf{R}(\boldsymbol{\beta})}
		{\partial\boldsymbol{\beta}}
	\end{equation}
	denotes the Jacobian matrix of the residual vector. Since the residual
	is constructed from the neural network representation and its
	derivatives, the Jacobian can be evaluated efficiently using automatic
	differentiation.
	
	The iteration is terminated when a prescribed convergence criterion is
	satisfied, yielding the optimized output weights
	$\boldsymbol{\beta}^{(\nu)}$. The corresponding neural network solution is
	then obtained from \eqref{uthetaII+}.
	
	\subsubsection{Optimal neural network approximation}
	The SCRINN loss \eqref{lossscrinn} can construct approximately orthonormal neural basis functions independently of any downstream approximation or PDE-solving task. Consequently, minimizing $\mathcal{L}_{\mathrm{SCRINN}}$ alone does not necessarily ensure that the resulting basis functions are well suited for a particular target function or PDE solution. In particular, orthogonality and normalization do not incorporate problem-specific information, such as local singularities, boundary layers, or other localized features of the target solution, and therefore do not guarantee that the learned basis functions can effectively capture such features.
	In addition, the Gram matrix in \eqref{2.3liner3} is computed from finitely many sampled points. Thus, minimizing $\mathcal{L}_{\mathrm{SCRINN}}$ enforces orthogonality and normalization only in a discrete sense. Prolonged training may further overfit the sampled points, resulting in a small training loss without preserving the desired properties over the entire domain. These two considerations motivate a stopping criterion that monitors the effectiveness of the learned basis functions in downstream applications while mitigating overfitting.

	To evaluate the performance of the trained neural network approximation, we introduce two quantitative metrics: the relative $L^2$ error and an MSE-based error. The relative $L^2$ error is defined by
	\begin{equation}\label{el2}
		E_{L^2}(u_\theta)
		:=
		\frac{\left(\sum\limits_{i=1}^{K}
			\left|u(\mathbf{x}_i)-u_{\theta}(\mathbf{x}_i)\right|^2\right)^{1/2}}
		{\left(\sum\limits_{i=1}^{K}
			\left|u(\mathbf{x}_i)\right|^2\right)^{1/2}},
	\end{equation}
	where $u(\mathbf{x}_i)$ and $u_{\theta}(\mathbf{x}_i)$ denote the exact solution and the neural network approximation evaluated at the sample points $\mathbf{x}_i$, respectively, and $K$ is the total number of evaluation points.
	
	In addition to the relative $L^2$ error in \eqref{el2}, we introduce an MSE-based error $E_{\mathrm{MSE}}$ as a task-dependent criterion for selecting the stopping point of Stage I. Unlike the SCRINN loss in \eqref{lossscrinn}, $E_{\mathrm{MSE}}$ is not involved in the optimization process. Instead, it is evaluated along the training trajectory and used to identify the point at which the learned neural basis functions provide the best performance for the downstream approximation or PDE-solving task.
	
	For function approximation, the MSE-based error is defined as
	\begin{equation*}
		E_{\mathrm{MSE}}(u_\theta)
		:=
		\frac{1}{N_{\mathrm{res}}}
		\sum_{i=1}^{N_{\mathrm{res}}}
		\left|
		u(\mathbf{x}_{i}^{\mathrm{res}})
		-
		u_{\theta}(\mathbf{x}_{i}^{\mathrm{res}})
		\right|^2,
	\end{equation*}
	while for PDE solving, it is defined as
	\begin{align*}
		E_{\mathrm{MSE}}(u_\theta)
		:= {}&
		\frac{1}{N_{\mathrm{res}}}
		\sum_{i=1}^{N_{\mathrm{res}}}
		\left|
		f(\mathbf{x}_{i}^{\mathrm{res}})
		-
		\mathscr{A}u_{\theta}(\mathbf{x}_{i}^{\mathrm{res}})
		\right|^2
		+
		\frac{1}{N_{\mathrm{bcs}}}
		\sum_{i=1}^{N_{\mathrm{bcs}}}
		\left|
		g(\mathbf{x}_{i}^{\mathrm{bcs}})
		-
		\mathscr{B}u_{\theta}(\mathbf{x}_{i}^{\mathrm{bcs}})
		\right|^2
		\\
		&+
		\frac{1}{N_{\mathrm{ics}}}
		\sum_{i=1}^{N_{\mathrm{ics}}}
		\left|
		\mathbf{h}(\mathbf{x}_{i}^{\mathrm{ics}})
		-
		\mathscr{I}u_{\theta}(\mathbf{x}_{i}^{\mathrm{ics}})
		\right|^2.
	\end{align*}
	Here, $\mathbf{x}_{i}^{\mathrm{res}}$, $\mathbf{x}_{i}^{\mathrm{bcs}}$, and
	$\mathbf{x}_{i}^{\mathrm{ics}}$ denote the interior residual, boundary, and
	initial sample points, respectively, and
	$N_{\mathrm{res}}$, $N_{\mathrm{bcs}}$, and $N_{\mathrm{ics}}$ denote the
	corresponding numbers of sample points. For elliptic PDEs, the
	initial-condition term is omitted.
	
	In practice, we optimize the SCRINN loss in \eqref{lossscrinn} for a fixed
	number of epochs $\nu_{\max}$, where $\nu_{\max}$ is typically a small
	integer (e.g., $\nu_{\max}\leq 20$). Rather than using the final iterate
	directly, we select the network parameters from the training trajectory
	according to the task-dependent criterion above. By default, we choose the
	epoch at which $E_{\mathrm{MSE}}$ attains its minimum. For function
	approximation problems, the relative $L^2$ error in \eqref{el2} may
	alternatively be used for this selection.
	Both $E_{\mathrm{MSE}}$ and $E_{L^2}$ are evaluated by forward propagation,
	with automatic differentiation used when derivatives of $u_{\theta}$ are
	required. Thus, their evaluation introduces only a relatively small
	additional computational cost.
	
	% More precisely, let $\hat{\theta}^{(\nu)}$ denote the network parameters at the $\nu$-th epoch obtained by minimizing \eqref{lossscrinn}. 
	We define the
	selected epoch index by
	\begin{equation}\label{nustar}
		\nu^*
		=
		\arg\min_{0\leq\nu\leq\nu_{\max}}
		E_{\mathrm{MSE}}\!\left(\hat{\theta}^{(\nu)}\right).
	\end{equation} 
	Then, the optimal neural network approximation is given by
	\begin{equation}\label{optimalnn}
		u_{\theta}^* =  u_{\theta}^{(\nu^*)}  (\mathbf{x}).
	\end{equation}

	\begin{remark}
		Minimizing the loss functions in both RINN and SCRINN trains the corresponding neural basis functions to be approximately orthonormal. Owing to the sparse architecture of SCRINN, the initial neural basis functions exhibit stronger linear independence, which enables the orthonormal basis functions to be learned more efficiently and accurately. Numerical results show that the stopping indicators for SCRINN typically reach their minima within only a few training epochs, whereas those for RINN often require substantially more epochs.
	\end{remark}

	The proposed SCRINN algorithm is summarized in
	Algorithm~\ref{algscrinn}.
	
	\begin{algorithm}[htbp]
		\caption{SCRINN for function approximation and PDE solving.}
		\label{algscrinn}
		\begin{algorithmic}[1]
			
			\REQUIRE
			Maximum number of epochs $\nu_{\max}$, learning rate $\eta$,
			penalty parameter $\varepsilon=1$, and Stage-I training points
			$\{\mathbf{x}_k\}_{k=1}^{K}$.
			For Stage II, data points $(\mathbf{X},u(\mathbf{X}))$ for function
			approximation, or residual, boundary, and initial points
			$\mathbf{x}^{\mathrm{res}}$, $\mathbf{x}^{\mathrm{bcs}}$, and
			$\mathbf{x}^{\mathrm{ics}}$, together with the corresponding
			functions $f$, $g$, and $\mathbf{h}$ for PDE problems.
			
			\ENSURE Optimal approximation $u_\theta^*$.
			
			\STATE Randomly initialize $\hat{\theta}^{(0)}$ and compute
			$\boldsymbol{\beta}^{(0)}$.
			\STATE Construct $u_\theta^{(0)}=\bigl(\boldsymbol{\beta}^{(0)}\bigr)^\top
			\boldsymbol{\phi}
			\bigl(\hat{\theta}^{(0)},\mathbf{x}\bigr)$ and set
			$\nu=1$, $\nu^*=0$.
			\STATE Compute $E_{\mathrm{MSE}}^* = E_{\mathrm{MSE}}(u_\theta^{(0)})$.
			
			\WHILE{$\nu \leq \nu_{\max}$}
			
			\algcomment{\textbf{Stage I: Learning approximately orthonormal neural basis functions}}
			
			\STATE Construct the neural basis functions
			\[
			\phi_j\bigl(\hat{\theta}^{(\nu-1)},\mathbf{x}_k\bigr),
			\qquad
			j=1,\ldots,N_L,\quad k=1,\ldots,K.
			\]
			
			\STATE Assemble the Gram matrix $\mathbf{G}$ using the prescribed
			quadrature rule in \eqref{discreteH1} or \eqref{discreteGaussian}.
			
			\STATE Evaluate the SCRINN loss
			$\mathcal{L}_{\mathrm{SCRINN}}$ defined in \eqref{lossscrinn}.
			
			\STATE Update the trainable parameters:
			\[
			\hat{\theta}^{(\nu)}
			=
			\hat{\theta}^{(\nu-1)}
			-
			\eta
			\nabla_{\hat{\theta}^{(\nu-1)}}
			\mathcal{L}_{\mathrm{SCRINN}}.
			\]
			
			\algcomment{\textbf{Stage II: Computing the output weights}}
			
			\STATE Freeze $\hat{\theta}^{(\nu)}$.
			
			\IF{function approximation or a linear PDE}
			
			\STATE Assemble the linear system
			\[
			\mathbf{H}\bigl(\hat{\theta}^{(\nu)},\mathbf{X}\bigr)
			\boldsymbol{\beta}
			=
			\mathbf{S}
			\]
			according to \eqref{2.2liner1approximation} or \eqref{2.2liner14}.
			
			\STATE Compute
			\[
			\boldsymbol{\beta}^{(\nu)}
			=
			\arg\min_{\boldsymbol{\beta}}
			\left\|
			\mathbf{H}\bigl(\hat{\theta}^{(\nu)},\mathbf{X}\bigr)
			\boldsymbol{\beta}
			-\mathbf{S}
			\right\|_2^2,
			\]
			with the solution given by \eqref{lspsol}.
			
			IF{a nonlinear PDE}\ELSIF{a nonlinear PDE}
			
			\STATE Construct the nonlinear residual vector
			$\mathbf{R}(\boldsymbol{\beta})$ according to \eqref{residualeqn}.
			
			\STATE Compute
			\[
			\boldsymbol{\beta}^{(\nu)}
			=
			\arg\min_{\boldsymbol{\beta}}
			\left\|
			\mathbf{R}(\boldsymbol{\beta})
			\right\|_2^2
			\]
			using the TRF nonlinear least-squares solver.
			
			\ENDIF
			
			\STATE Construct
			\[
			u_{\theta}^{(\nu)}(\mathbf{x})
			=
			\bigl(\boldsymbol{\beta}^{(\nu)}\bigr)^\top
			\boldsymbol{\phi}
			\bigl(\hat{\theta}^{(\nu)},\mathbf{x}\bigr).
			\]
			
			\STATE Evaluate $E_{\mathrm{MSE}}(u_\theta^{(\nu)})$.
			
			\IF{$E_{\mathrm{MSE}}(u_\theta^{(\nu)})
				<
				E_{\mathrm{MSE}}^*$}
			\STATE Set $\nu^*=\nu$,
			$u_\theta^*=u_\theta^{(\nu)}$, and $E_{\mathrm{MSE}}^* = E_{\mathrm{MSE}}(u_\theta^{(\nu)})$.
			\ENDIF
			
			\STATE Set $\nu=\nu+1$.
			
			\ENDWHILE
			
			\RETURN $u_\theta^*$.
			
		\end{algorithmic}
	\end{algorithm}

	\begin{remark}\label{RINN+}
		To further assess the effectiveness of the proposed loss function \eqref{lossscrinn}, we investigate its applicability beyond the SCRINN framework. Specifically, we incorporate the proposed loss formulation into the existing RINN-es model while keeping its original network architecture unchanged. The resulting model is referred to as \textbf{RINN+}.
		
		RINN+ has exactly the same network architecture, connectivity pattern, and number of trainable parameters as the original RINN-es; the only difference is the loss function used for training. This controlled modification provides a direct and fair comparison between RINN-es and RINN+, allowing the effect of the proposed loss formulation to be evaluated independently of architectural changes.
	\end{remark}

	\section{Numerical examples}\label{sec:examples}
	
	In this section, we present numerical experiments to assess the performance of the proposed SCRINN for function approximation and the numerical solution of PDEs. We compare SCRINN with the RINN-es method proposed in \cite{PengHuangYi2026} and its enhanced variant, RINN+, introduced in \Cref{RINN+}. All experiments are conducted on an NVIDIA A100-SXM4-80GB GPU for both training and evaluation, except for \Cref{ex:2ha}, which involves nonlinear least-squares optimization. For this example, the nonlinear least-squares problem in Stage II is solved on a CPU (Intel Xeon Platinum 8358 @ 2.60 GHz) using the SciPy library.

	Unless otherwise specified, all \(m\)-dimensional SCRINN experiments employ the default architecture
	\[
	N_{nn}=[m,4000,4000,1],
	\]
	where the weights and biases are initialized from the uniform distribution \(\mathcal{U}(-1,1)\), and the sine function $sin(\cdot)$ is adopted as the activation function.
	In Stage I, we set \(\varepsilon=1\) in the SCRINN loss \eqref{lossscrinn} and train the network using the Adam optimizer with a learning rate of \(\eta=5\times10^{-2}\) for at most \(\nu_{\max}=20\) epochs. For the two-dimensional experiments, the computational domain \(\mathcal{D}\) is partitioned into \(16\times16\) elements with \(4\times4\) Gaussian quadrature points per element, yielding \(4096\) interior collocation points. For the three-dimensional experiments, \(\mathcal{D}\) is partitioned into \(4\times4\times4\) elements with \(5\times5\times5\) Gaussian quadrature points per element, also resulting in \(8000\) interior collocation points.
	In Stage II, an additional \(6144\) uniformly distributed points are sampled on the boundary \(\partial\mathcal{D}\) for the least-squares solve in both the two- and three-dimensional experiments.
	Finally, the trained SCRINN is evaluated on uniformly distributed test points in \(\overline{\mathcal{D}}\). Specifically, \(40401\) test points are used in the two-dimensional experiments, whereas \(125000\) test points are used in the three-dimensional experiments. The relative \(L^2\) error \(E_{L^2}\) is computed according to \eqref{el2}.
	
	\subsection{SCRINN for Function Approximation}
	We first investigate the function approximation capability of SCRINN and assess how its sparse connectivity architecture and orthogonalization strategy contribute to the representational capacity and training stability of neural networks.

	\begin{example} \label{ex:Drop-Wave}
		The Drop-Wave function is a highly oscillatory function with a sharp global maximum at the origin and concentric oscillatory structures away from the origin. It is defined by
		\begin{equation} \label{3.1.2liner1}
			u\left( x,y \right)
			=
			\frac{1+\cos \left( 12\sqrt{x^2+y^2} \right)}
			{0.5\left( x^2+y^2 \right) +2}.
		\end{equation}
		We consider its approximation on the square domain $\Omega=[-2,2]^2$.
	\end{example}
	
\begin{figure}[htbp]
		\centering
		\begin{subfigure}[b]{0.36\linewidth}
			\centering
			\includegraphics[width=0.9\linewidth]{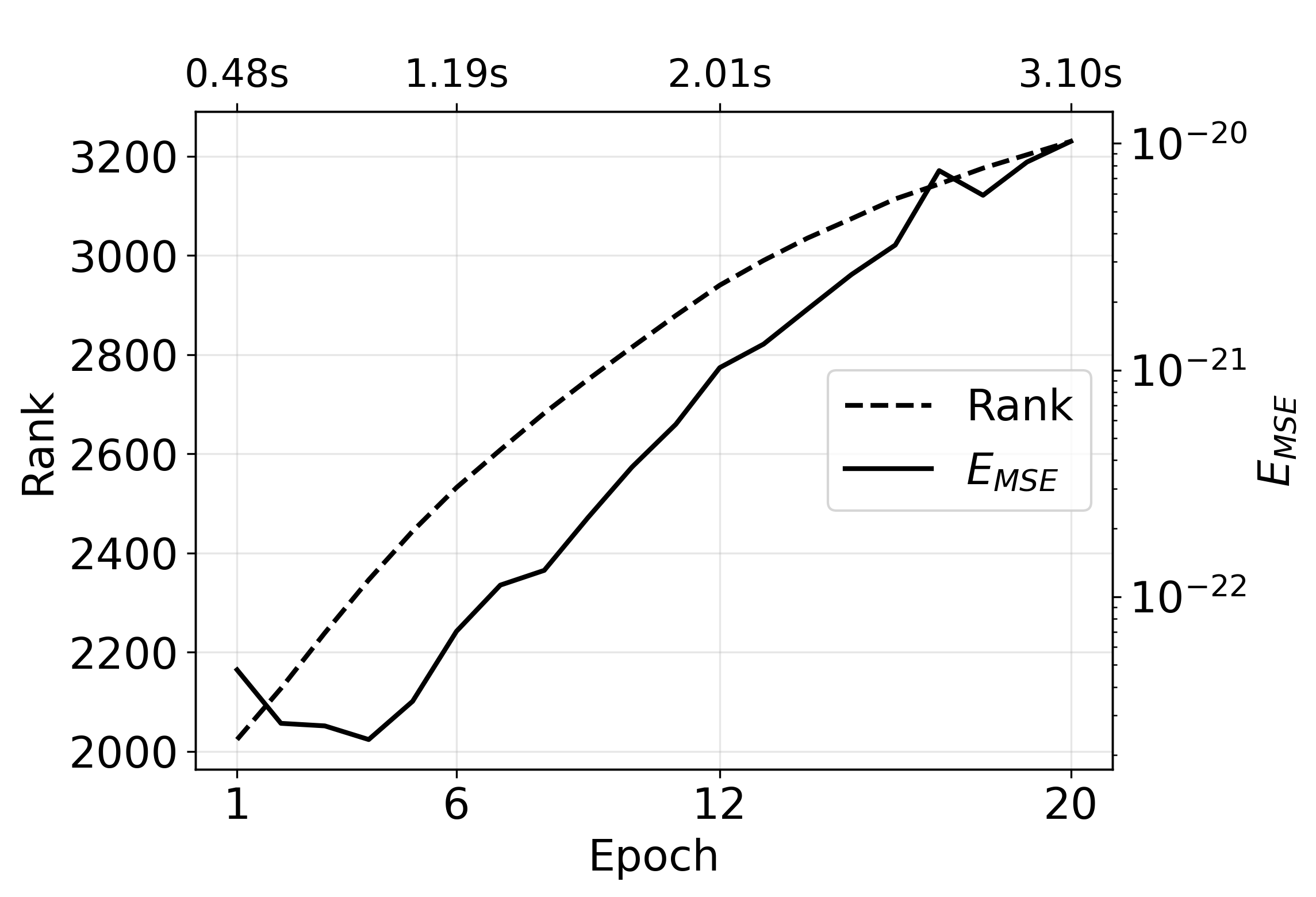}
			\caption{}
			\label{fig:Drop-Wave-Training}
		\end{subfigure}
		\hfill
		\begin{subfigure}[b]{0.32\linewidth}
			\centering
			\includegraphics[width=0.9\linewidth]{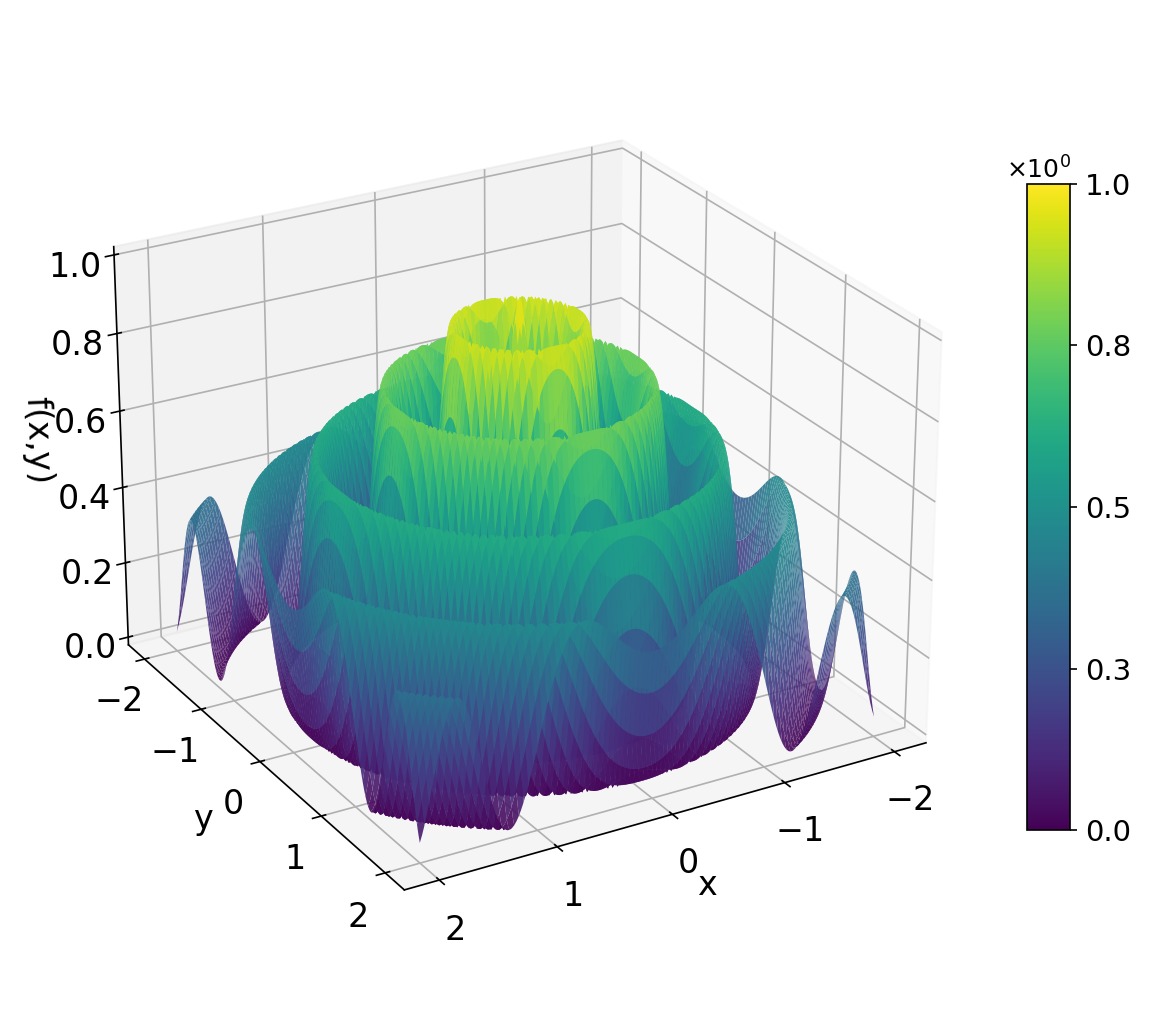}
			\caption{}
			\label{fig:Drop-Wave-Approximation}
		\end{subfigure}
		\hfill
		\begin{subfigure}[b]{0.30\linewidth}
			\centering
			\includegraphics[width=0.9\linewidth]{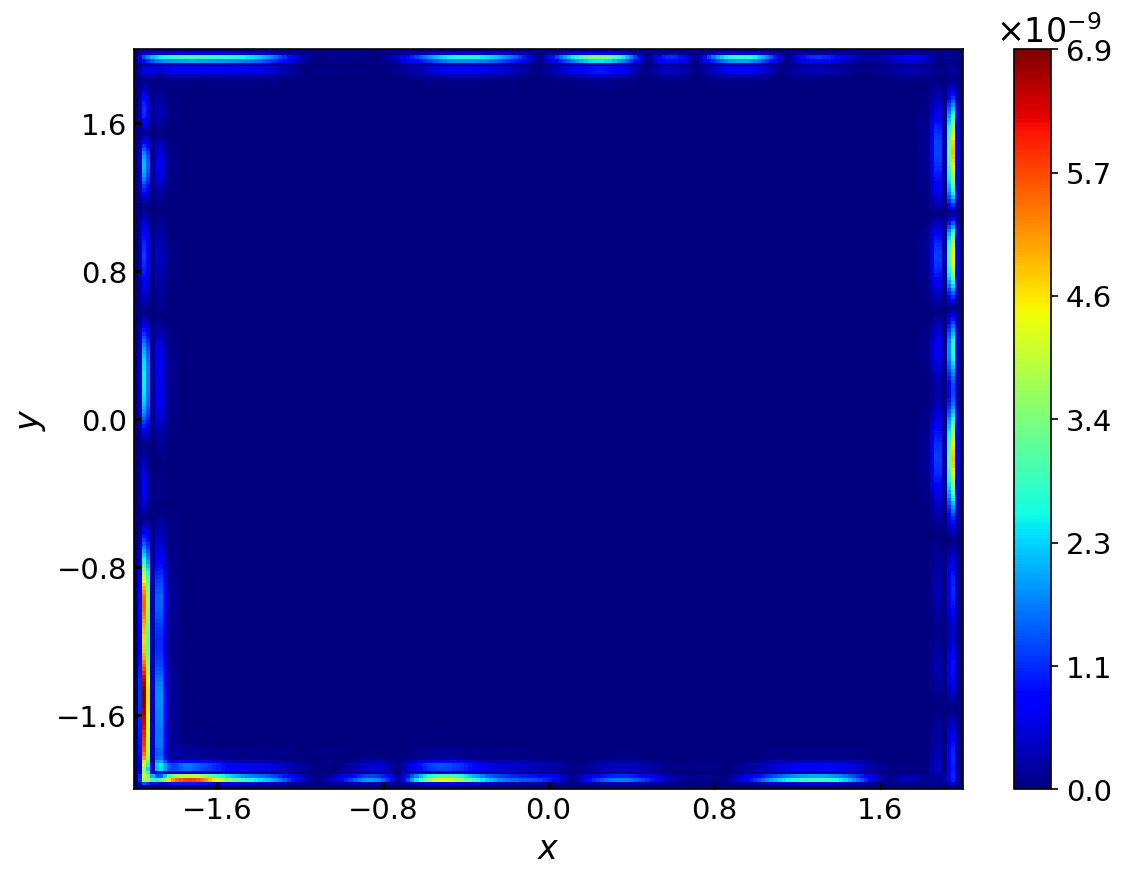}
			\caption{}
			\label{fig:Drop-Wave-Error}
		\end{subfigure}
		
		\caption{\Cref{ex:Drop-Wave}. Training and approximation results of SCRINN for the Drop-Wave function.
			(a) Evolution of $E_{\mathrm{MSE}}$ and the effective rank of the hidden-layer output matrix;
			(b) SCRINN approximation;
			(c) pointwise absolute error.}
		\label{fig:DropWaveComparison}
	\end{figure}

	\Cref{fig:Drop-Wave-Training} shows the evolution of the hidden-layer output matrix rank and the MSE-based error $E_{\mathrm{MSE}}$ during Stage I training of SCRINN. The rank increases progressively during training, indicating that the hidden-layer output space is gradually enriched as the training proceeds. At the same time, $E_{\mathrm{MSE}}$ remains at a low level, on the order of $10^{-20}$, with only small fluctuations. These observations indicate that the rank enrichment is achieved while maintaining stable training behavior.
	
	\Cref{fig:Drop-Wave-Approximation} and \Cref{fig:Drop-Wave-Error} present the SCRINN approximation and the corresponding pointwise absolute error, respectively. The approximation captures the main features of the target function, including the sharp maximum at the origin and the surrounding concentric oscillatory structures. The pointwise absolute error remains small throughout the computational domain, indicating that SCRINN provides an accurate approximation of this highly oscillatory function.
	
	\Cref{tab:3.1.2tab1} compares the approximation performance of ELM, RINN-es, and SCRINN for different numbers of neural basis functions. For SCRINN, increasing the number of basis functions from 2000 to 6000 reduces the relative $L^2$ error from $1.3524\times10^{-6}$ to $1.0479\times10^{-10}$. Over the same range, the effective rank of the hidden-layer output matrix increases from 1030 to 3444. This behavior is consistent with the role of the sparse connection architecture in progressively enriching the neural basis.
	
	In contrast, ELM does not provide a satisfactory approximation for this test problem, with relative $L^2$ errors remaining at the level of $10^1$--$10^2$ as the number of basis functions increases. RINN-es achieves relative $L^2$ errors of $4.7098\times10^{-8}$ and $2.8560\times10^{-8}$ with 2000 and 4000 basis functions, respectively, but does not maintain this accuracy for 6000 basis functions, for which the error increases to $8.5361\times10^{-3}$. In addition to the approximation accuracy, SCRINN requires substantially less training time. For 6000 basis functions, the training time is reduced from 61.9 seconds for RINN-es to 4.8 seconds for SCRINN. These results demonstrate that SCRINN can achieve high approximation accuracy while maintaining a progressively increasing effective rank and substantially reducing the training cost.
	
	\Cref{tab:nihelossdiif} examines the sensitivity of RINN-es and SCRINN to the hyperparameter $\varepsilon$ under the original loss \eqref{lossf} and the proposed loss \eqref{lossscrinn}. For RINN-es, the original loss exhibits a strong dependence on $\varepsilon$. Its relative $L^2$ error is $2.8560\times10^{-8}$ for $\varepsilon=0.01$, but increases to $1.1729\times10^{-2}$ and $2.2113\times10^{-2}$ for $\varepsilon=0.5$ and $1.0$, respectively. With the proposed loss, the corresponding RINN+ results remain close to $2\times10^{-8}$ for all three values of $\varepsilon$. Thus, the proposed loss substantially reduces the sensitivity of the RINN formulation to $\varepsilon$.
	
	For SCRINN, both loss functions yield relative $L^2$ errors of the order of $10^{-9}$ for all tested values of $\varepsilon$. The proposed loss gives errors of $1.1415\times10^{-9}$, $1.1976\times10^{-9}$, and $1.1133\times10^{-9}$ for $\varepsilon=0.01$, $0.5$, and $1.0$, respectively. In particular, the smallest error is obtained with $\varepsilon=1.0$. These results indicate that SCRINN is less sensitive to $\varepsilon$ than RINN-es for this test problem, while the proposed loss provides consistent approximation accuracy across the tested parameter range. Based on these results, we set $\varepsilon=1.0$ in the subsequent experiments.

	\begin{table}[htbp]
		\centering
		\caption{\Cref{ex:Drop-Wave}. Approximation results of different models for the Drop-Wave function.}
		\label{tab:3.1.2tab1}
		\begin{tabular}{ccccc}
			\toprule
			Method & Architecture & $E_{L^2}$ & Rank & Training time(s) \\
			\midrule
			
			\multirow{3}{*}{ELM}
			& $[2,512,2000,1]$ & $3.9517\times 10^{+1}$ & 2000 & 3.9 \\
			& $[2,512,4000,1]$ & $1.1090 \times 10^{+2}$ & 4000 & 19.9 \\
			& $[2,512,6000,1]$ &  $5.7723\times 10^{+1}$ & 4096 & 20.2 \\
			
			\midrule
			
			\multirow{3}{*}{RINN-es}
			& $[2,512,2000,1]$ & $4.7098\times 10^{-8}$ & 1448 & 36.1 \\
			& $[2,512,4000,1]$ & $2.8560 \times 10^{-8}$ & 2419 & 57 \\
			& $[2,512,6000,1]$ & $8.5361 \times 10^{-3}$ & 4096 & 61.9 \\
			
			\midrule
			
			\multirow{3}{*}{SCRINN}
			& $[2,2000,2000,1]$ & $1.3524 \times 10^{-6}$ & 1030 & 1.4 \\
			& $[2,4000,4000,1]$ & $1.1133 \times 10^{-9}$ & 2346 & 3.1 \\
			& $[2,6000,6000,1]$ & $1.0479 \times 10^{-10}$ & 3444 & 4.8 \\
			%& $[2,8000,8000,1]$ & $1.6625 \times 10^{-8}$ & 4096 & 5.3 \\
			\bottomrule
		\end{tabular}
	\end{table}
	
	\begin{table}[htpb]
		\centering
		\caption{\Cref{ex:Drop-Wave}. Influence of different loss functions and $\varepsilon$ values on $E_{L^2}$ of RINN-es and SCRINN for the Drop-Wave function. Boldface entries correspond to the results obtained by RINN+.}
		\label{tab:nihelossdiif}
		\begin{tabular}{ccccc}
			\toprule
			\text{Method} & \text{Architecture} &
			\boldmath{$\varepsilon$} &
			Original Loss \eqref{lossf} &
			Proposed Loss \eqref{lossscrinn} \\
			\midrule
			
			\multirow{3}{*}{RINN-es} 
			& \multirow{3}{*}{$[2,512,4000,1]$}
			& 0.01 & $2.8560\times10^{-8}$ & $\mathbf{1.9955\times10^{-8}}$ \\
			& & 0.5 & $1.1729\times10^{-2}$ & $\mathbf{2.0387\times10^{-8}}$ \\
			& & 1.0 & $2.2113\times10^{-2}$ & $\mathbf{1.8674\times10^{-8}}$ \\
			
			\midrule
			
			\multirow{3}{*}{SCRINN} 
			& \multirow{3}{*}{$[2,4000,4000,1]$}
			& 0.01 & $6.2786\times10^{-10}$ & $1.1415\times10^{-9}$ \\
			& & 0.5 & $1.4230\times10^{-9}$ & $1.1976\times10^{-9}$ \\
			& & 1.0 & $1.8848\times10^{-9}$ & $1.1133\times10^{-9}$ \\
			
			\bottomrule
		\end{tabular}
	\end{table}

	\begin{example}\label{ex:Michalewicz}
		The Michalewicz function is a classical multimodal benchmark function commonly used to assess numerical methods for nonlinear approximation and optimization. Its complexity is controlled by the parameter $n$, which determines the sharpness of its local minima and the degree of oscillation. As $n$ increases, the function develops increasingly localized and oscillatory structures, providing a challenging test for the approximation of nonlinear function landscapes. It is defined by
		\begin{equation}
			\label{3.1.4liner1}
			f(\mathbf{x})
			=
			\sum_{i=1}^{d}
			\sin(x_i)
			\left[
			\sin\left(\frac{i x_i^2}{\pi}\right)
			\right]^{2n},
		\end{equation}
		where $\mathbf{x}=(x_1,\ldots,x_d)$ denotes the input vector, $d$ is the dimension of the function, and $n$ is a positive parameter controlling the complexity of the function. We consider the approximation of the two-dimensional Michalewicz function on the domain $\Omega=[1,3]^2$ with $n=10$.
	\end{example}
	
\begin{figure}[htbp]
		\centering
		\begin{subfigure}[b]{0.36\linewidth}
			\centering
			\includegraphics[width=0.9\linewidth]{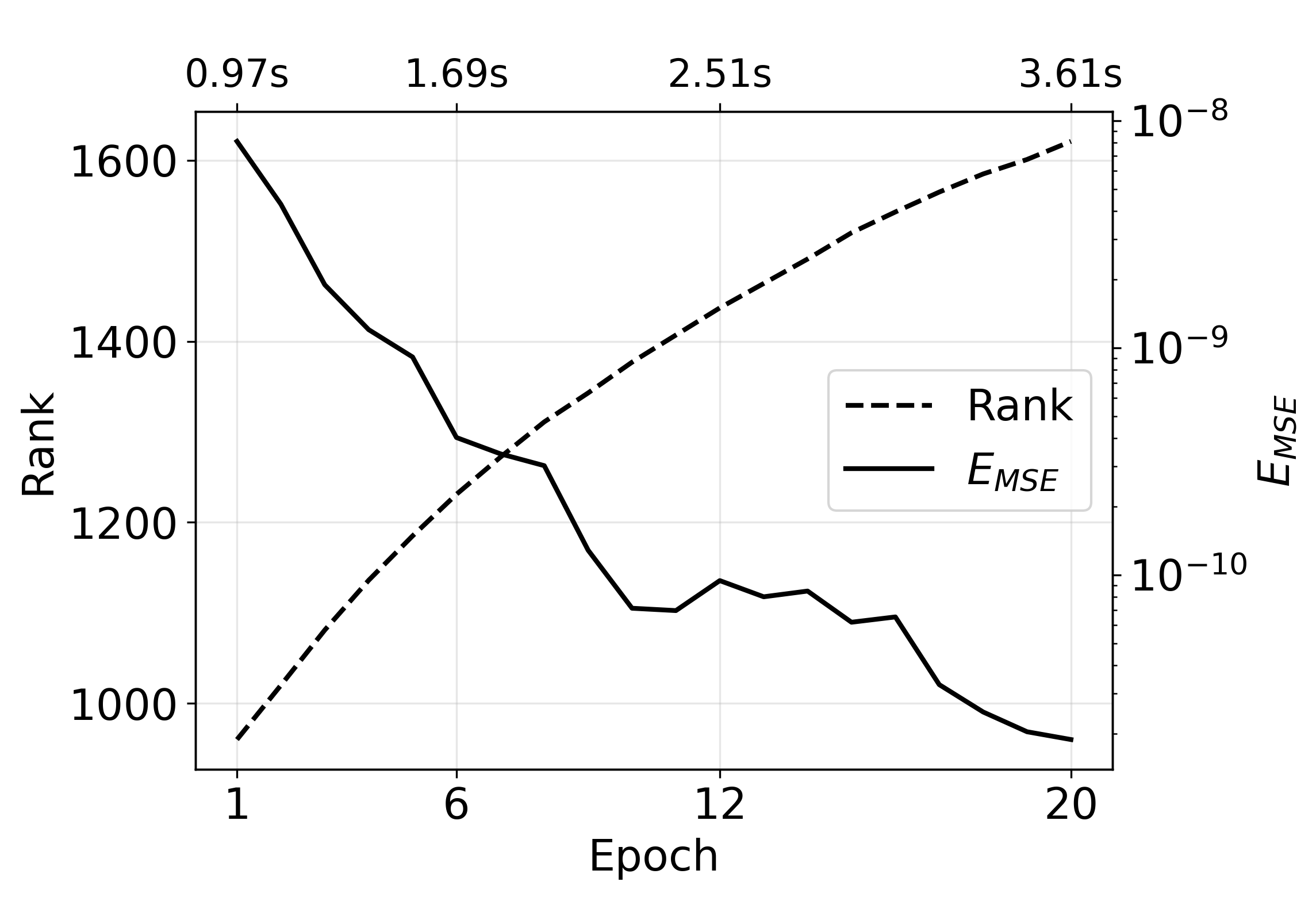}
			\caption{}
			\label{fig:Michalewiczjieguom=10quxian}
		\end{subfigure}
		\hfill
		\begin{subfigure}[b]{0.32\linewidth}
			\centering
			\includegraphics[width=0.9\linewidth]{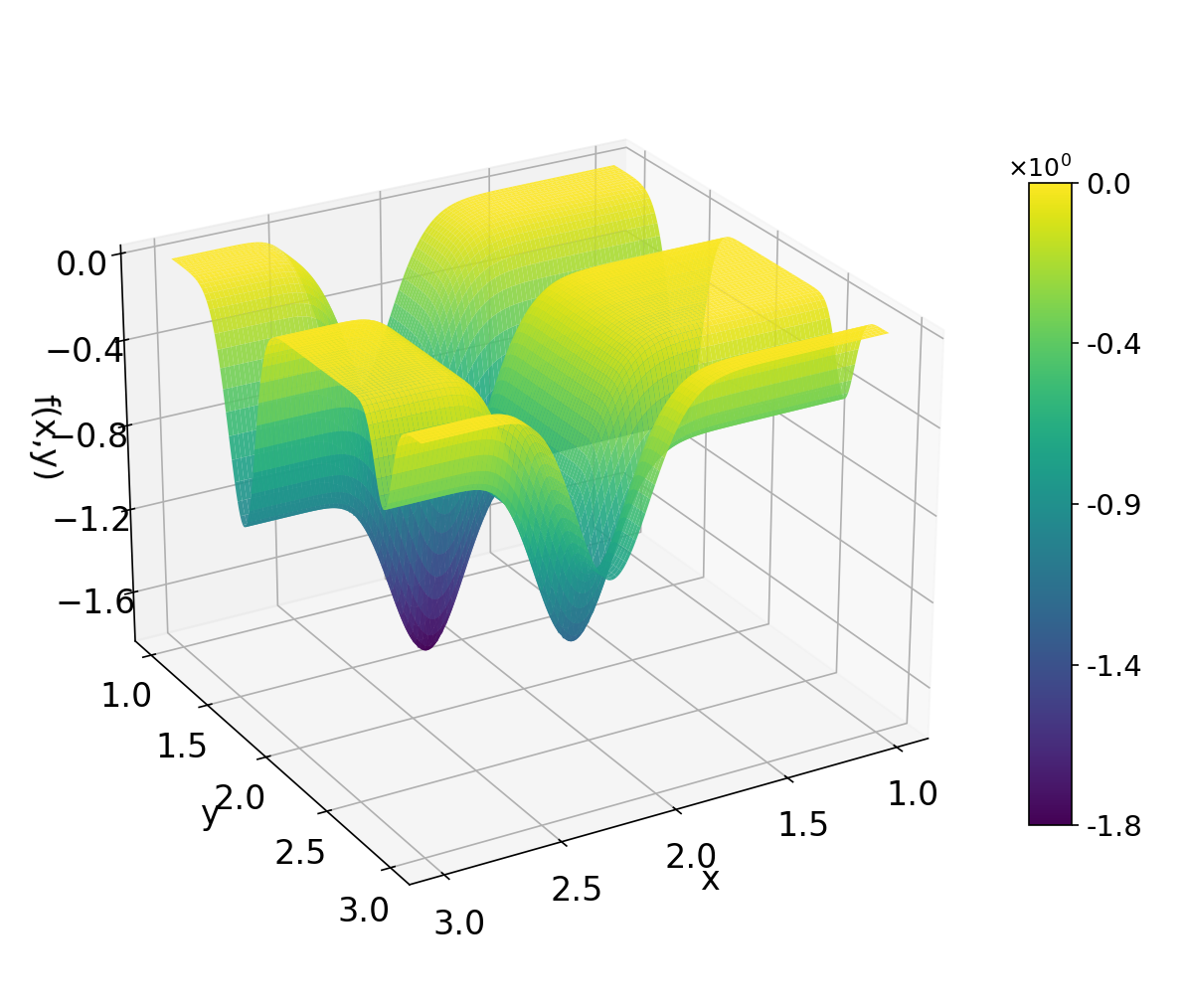}
			\caption{}
			\label{fig:Michalewiczm=10}
		\end{subfigure}
		\hfill
		\begin{subfigure}[b]{0.30\linewidth}
			\centering
			\includegraphics[width=0.9\linewidth]{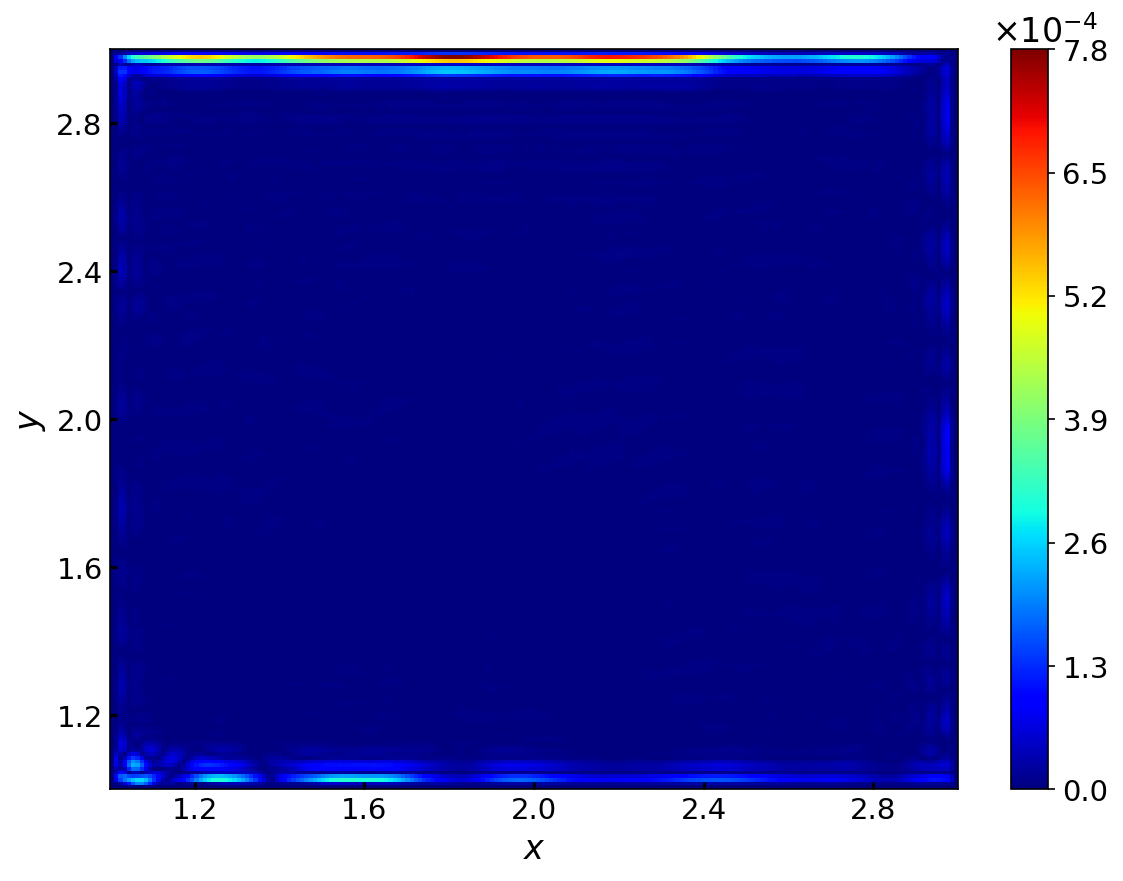}
			\caption{}
			\label{fig:Michalewiczjieguom=10}
		\end{subfigure}
		
		\caption{\Cref{ex:Michalewicz}. Training and approximation results of SCRINN for the Michalewicz function.
			(a) Evolution of $E_{\mathrm{MSE}}$ and the effective rank of the hidden-layer output matrix;
			(b) SCRINN approximation;
			(c) pointwise absolute error.}
		\label{fig:Michalewiczcomparison}
	\end{figure}

	\Cref{fig:Michalewiczjieguom=10quxian} shows the evolution of the MSE-based error $E_{\mathrm{MSE}}$ and the rank of the hidden-layer output matrix during Stage I training of SCRINN with 4000 neural basis functions. The effective rank increases rapidly during the initial stage of training and continues to grow as the training proceeds, indicating a progressive enrichment of the hidden-layer output space. Meanwhile, $E_{\mathrm{MSE}}$ decreases by several orders of magnitude, reflecting the improvement in the approximation of the neural basis functions. In the later stage of training, the increase in rank becomes more gradual and the decrease in $E_{\mathrm{MSE}}$ slows down. This behavior suggests that the training progressively enriches the neural basis while reducing the approximation error.
	
	\Cref{fig:Michalewiczm=10} presents the SCRINN approximation of the Michalewicz function. The approximation captures the main features of the target function, including its nonlinear and localized oscillatory structures. The corresponding pointwise absolute error in \Cref{fig:Michalewiczjieguom=10} remains small over the computational domain. The results show that SCRINN can represent the nonlinear and multimodal structure of the Michalewicz function while maintaining stable training behavior.
	
	\Cref{tab:3.1.4tab1} summarizes the approximation results of ELM, RINN-es, and SCRINN for different numbers of neural basis functions. For SCRINN, increasing the number of basis functions from 2000 to 6000 results in a progressive increase in the effective rank of the hidden-layer output matrix, from 890 to 2933. Meanwhile, the relative $L^2$ error decreases from $1.7179\times10^{-3}$ to $5.9261\times10^{-5}$. The decrease in the approximation error together with the increase in effective rank indicates that the additional neural basis functions provide further representational capacity for this test problem.
	
	In contrast, ELM and RINN-es yield relative $L^2$ errors at the level of $10^{1}$--$10^{2}$ for the tested network sizes. In comparison, SCRINN achieves substantially smaller approximation errors while requiring less training time. For 6000 neural basis functions, the training time of SCRINN is 4.4 seconds, compared with 7.7 and 19.6 seconds for ELM and RINN-es, respectively. These results indicate that the sparse connection architecture enables SCRINN to obtain a more effective neural representation with relatively low computational cost.

	\begin{example}\label{ex:Rosenbrock}
		We next consider the Rosenbrock function
		\begin{equation}
			\label{3.1.3liner1}
			u(x,y)=\left(1-x\right)^2+100\left(y-x^2\right)^2,
		\end{equation}
		on $\Omega=[-6,6]^2$. The function provides a nonseparable test problem
		with strong nonlinear coupling and anisotropic variation associated with
		its narrow curved valley.
	\end{example}
	
\begin{figure}[htbp]
		\centering
		\begin{subfigure}[b]{0.36\linewidth}
			\centering
			\includegraphics[width=0.9\linewidth]{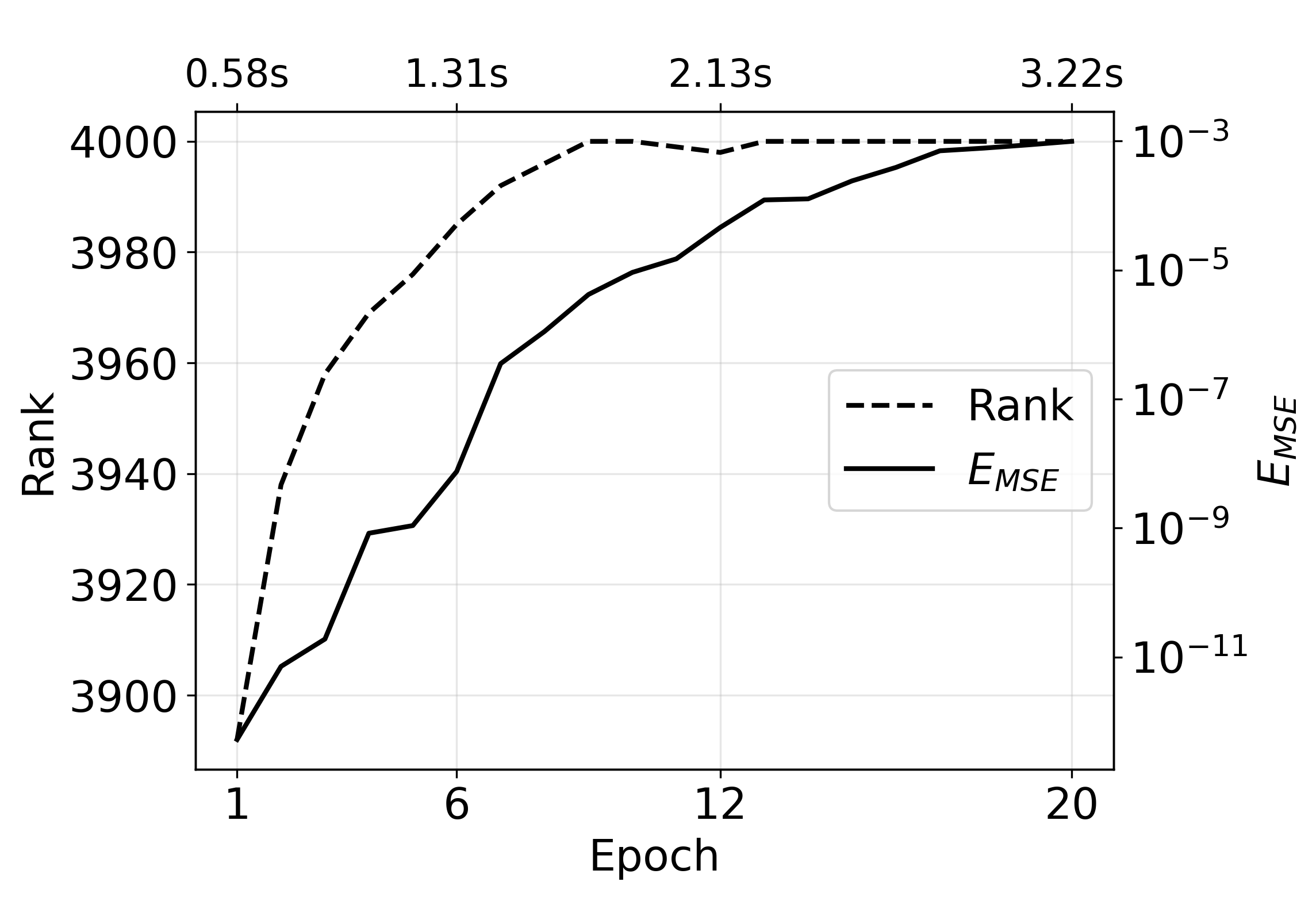}
			\caption{}
			\label{fig:Rosenbrockquxian}
		\end{subfigure}
		\hfill
		\begin{subfigure}[b]{0.32\linewidth}
			\centering
			\includegraphics[width=0.9\linewidth]{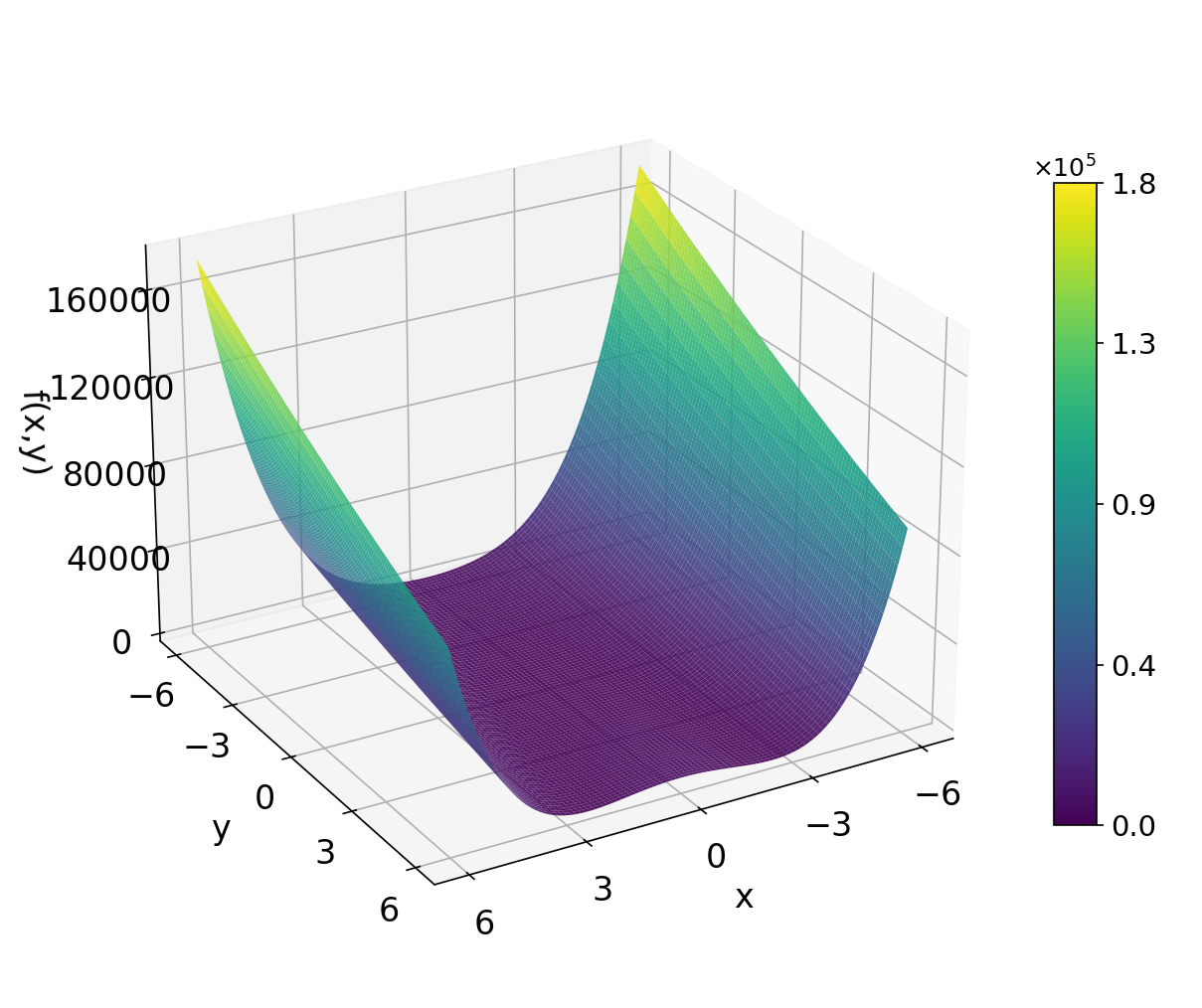}
			\caption{}
			\label{fig:Rosenbrock}
		\end{subfigure}
		\hfill
		\begin{subfigure}[b]{0.3\linewidth}
			\centering
			\includegraphics[width=0.9\linewidth]{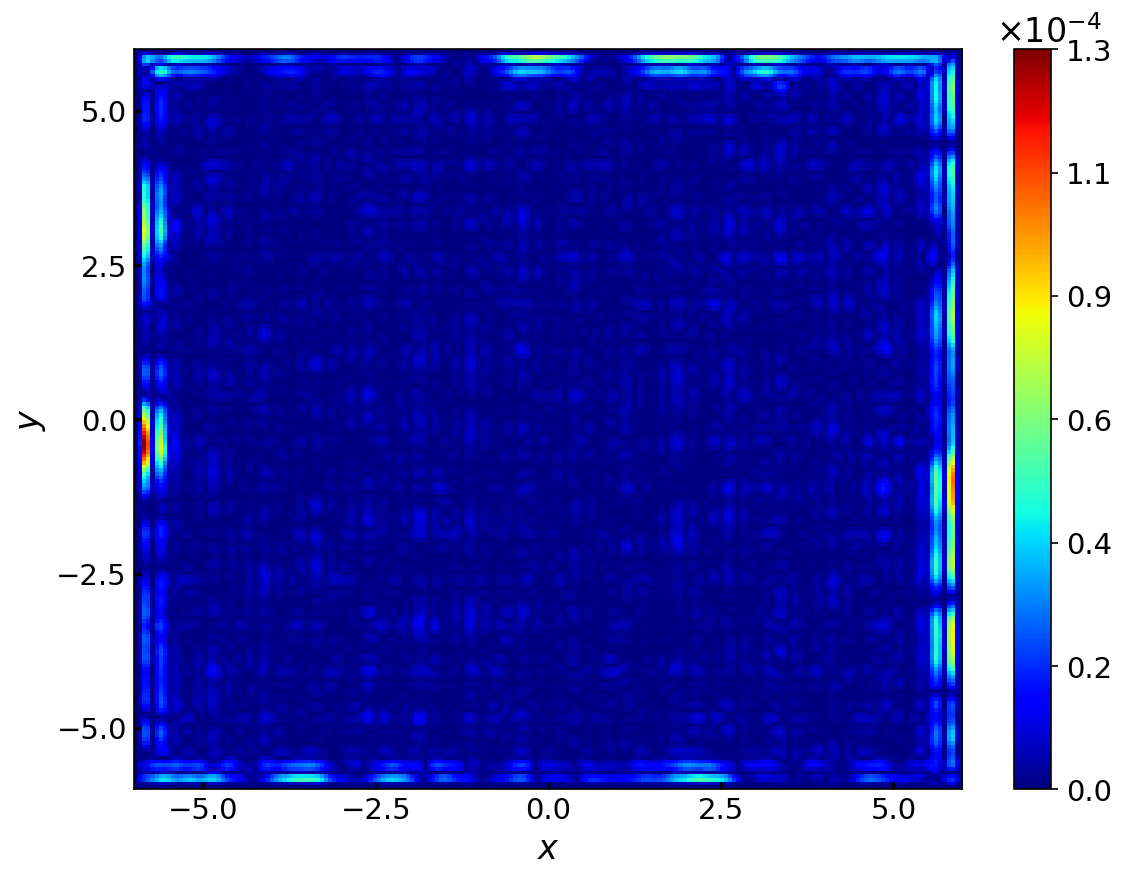}
			\caption{}
			\label{fig:Rosenbrockjieguo}
		\end{subfigure}
		\caption{\Cref{ex:Rosenbrock}. Training process and approximation results of SCRINN on the Rosenbrock function.
			(a) Evolution of $E_{\mathrm{MSE}}$ and the rank of the hidden-layer output matrix;
			(b) SCRINN approximation; 
			(c) pointwise absolute error.}
		\label{fig:Rosenbroc}
	\end{figure}

	\Cref{fig:Rosenbrockquxian} shows the evolution of the approximation error
	$E_{\mathrm{MSE}}$ and the rank of the hidden-layer output matrix during
	Stage I training with 4000 neural basis functions. The initial neural basis
	already has a high rank, and the rank rapidly approaches its maximum during
	the early training iterations. Meanwhile, the approximation error increases
	during this initial stage. This behavior indicates that, for this test
	problem, the randomly initialized basis already provides a sufficiently rich
	representation, while further Stage I optimization primarily changes the
	basis structure rather than improving the approximation. 
	% We therefore use
	% only a limited number of Stage I iterations for the subsequent experiments.
	
	\Cref{fig:Rosenbrock} and \Cref{fig:Rosenbrockjieguo} show the SCRINN
	approximation and the corresponding pointwise absolute error obtained with
	4000 neural basis functions. The approximation reproduces the main structure
	of the Rosenbrock function, including its narrow curved valley, and the
	pointwise error remains small over the computational domain.
	
	\Cref{tab:3.1.3tab1} compares SCRINN with ELM and RINN-es for different
	numbers of neural basis functions. For SCRINN, the relative $L^2$ error
	decreases from $1.7898\times10^{-5}$ to $1.8341\times10^{-10}$ as the
	number of basis functions increases from 2000 to 6000, while the effective
	rank increases from 2000 to 3991. ELM yields relative $L^2$ errors close to
	$1$ for all tested network sizes, whereas RINN-es achieves errors on the
	order of $10^{-2}$. SCRINN also requires substantially less training time
	than RINN-es. Thus, for this test problem, the sparse architecture provides
	a more accurate approximation with a lower computational cost than the
	fully connected RINN-es architecture.

	\begin{table}[htbp]
		\centering
		\caption{\Cref{ex:Rosenbrock}. Approximation results of different models for the Rosenbrock function.}
		\label{tab:3.1.3tab1}
		\begin{tabular}{ccccc}
			\toprule
			Method & Architecture & $E_{L^2}$ & Rank & Time(s) \\
			\midrule
			
			\multirow{3}{*}{ELM}
			& $[2,512,2000,1]$ & $9.9993 \times 10^{-1}$ & 2000 & 3.4 \\
			& $[2,512,4000,1]$ & $9.9960 \times 10^{-1}$ & 4000 & 7.2 \\
			& $[2,512,6000,1]$ & $1.0000 \times 10^{0}$ & 4096 & 7.6 \\
			
			\midrule
			
			\multirow{3}{*}{RINN-es}
			& $[2,512,2000,1]$ & $1.9245 \times 10^{-2}$ & 2000 & 12.7 \\
			& $[2,512,4000,1]$ & $6.2952 \times 10^{-2}$ & 4000 & 15.5 \\
			& $[2,512,6000,1]$ & $5.0464\times 10^{-2}$ & 4096 & 27.1 \\
			
			\midrule
			
			\multirow{3}{*}{SCRINN}
			& $[2,2000,2000,1]$ & $1.7898 \times 10^{-5}$ & 2000 & 1.4 \\
			& $[2,4000,4000,1]$ & $2.0373 \times 10^{-10}$ & 3892 & 3.2 \\
			& $[2,6000,6000,1]$ & $1.8341 \times 10^{-10}$ & 3991 & 4.9 \\
			%& $[2,8000,8000,1]$ & $6.6860 \times 10^{-5}$ & 4096 & 7.12 \\
			\bottomrule
		\end{tabular}
	\end{table}

	\begin{example}\label{ex:3Dfunc2}
		We further examine the approximation capability of SCRINN for the three-dimensional function defined on $\Omega=[-1,1]^3$,
		\begin{equation}
			\label{3dlaxs5}
			u(\mathbf{x})
			=
			1-
			\cos(2\pi x_1)
			\cos(2\pi(x_2-x_1))
			\cos(2\pi(x_3-x_2))
			e^{-0.1(x_1^2+x_2^2+x_3^2)}.
		\end{equation}
		The function combines oscillatory behavior, successive coupling of neighboring variables, and a Gaussian envelope, providing a nonseparable three-dimensional approximation problem.
	\end{example}
	
	\Cref{tab:3dtab1} reports the relative $L^2$ error, rank, and training
	time for ELM, RINN-es, and SCRINN with different numbers of neural basis
	functions. For SCRINN, the relative $L^2$ error decreases from
	$4.7484\times10^{-2}$ to $1.4124\times10^{-3}$ and
	$3.6651\times10^{-4}$ as the number of basis functions increases from
	2000 to 4000 and 6000, respectively. The corresponding ranks
	increase from 2000 to 3992 and 5891. RINN-es also exhibits decreasing errors
	with increasing basis size, but its errors remain larger than those of SCRINN
	at all three network sizes. In particular, for 4000 and 6000 basis functions,
	SCRINN achieves relative $L^2$ errors of $1.4124\times10^{-3}$ and
	$3.6651\times10^{-4}$, compared with $1.3985\times10^{-2}$ and
	$3.7126\times10^{-3}$ for RINN-es.
	
	The computational cost shows a similar advantage for SCRINN. Its training
	time increases from $10.6$ to $18.2$ seconds as the number of basis functions
	increases from 2000 to 6000, whereas the corresponding times for RINN-es
	increase from $11.8$ to $126.3$ seconds. ELM is faster but produces
	substantially larger errors, which increase with the network size. These
	results indicate that the sparse connectivity of SCRINN allows the use of
	larger neural bases while maintaining relatively low training costs.
	
	\begin{table}[htbp]
		\centering
		\caption{\Cref{ex:3Dfunc2}. Relative $L^2$ error, effective rank, and training time for the three-dimensional function.}
		\label{tab:3dtab1}
		\begin{tabular}{ccccc}
			\toprule
			Method & Architecture & $E_{L^2}$ & Effective rank & Time (s) \\
			\midrule
			
			\multirow{3}{*}{ELM}
			& $[3,512,2000,1]$
			& $1.6357\times10^{1}$ & 2000 & 4.7\\
			& $[3,512,4000,1]$
			& $2.8711\times10^{1}$ & 4000 & 8.0\\
			& $[3,512,6000,1]$
			& $5.0550\times10^{1}$ & 6000 & 14.2\\
			
			\cmidrule{1-5}
			
			\multirow{3}{*}{RINN-es}
			& $[3,512,2000,1]$
			& $8.7748\times10^{-2}$ & 2000 & 11.8\\
			& $[3,512,4000,1]$
			& $1.3985\times10^{-2}$ & 4000 & 59.1\\
			& $[3,512,6000,1]$
			& $3.7126\times10^{-3}$ & 5696 & 126.3\\
			
			\cmidrule{1-5}
			
			\multirow{3}{*}{SCRINN}
			& $[3,2000,2000,1]$
			& $4.7484\times10^{-2}$ & 2000 & 10.6\\
			& $[3,4000,4000,1]$
			& $1.4124\times10^{-3}$ & 3992 & 12.6\\
			& $[3,6000,6000,1]$
			& $3.6651\times10^{-4}$ & 5891 & 18.2\\
			
			\bottomrule
		\end{tabular}
	\end{table}

\begin{figure}[htbp]
		\centering
		\includegraphics[width=0.85\linewidth]{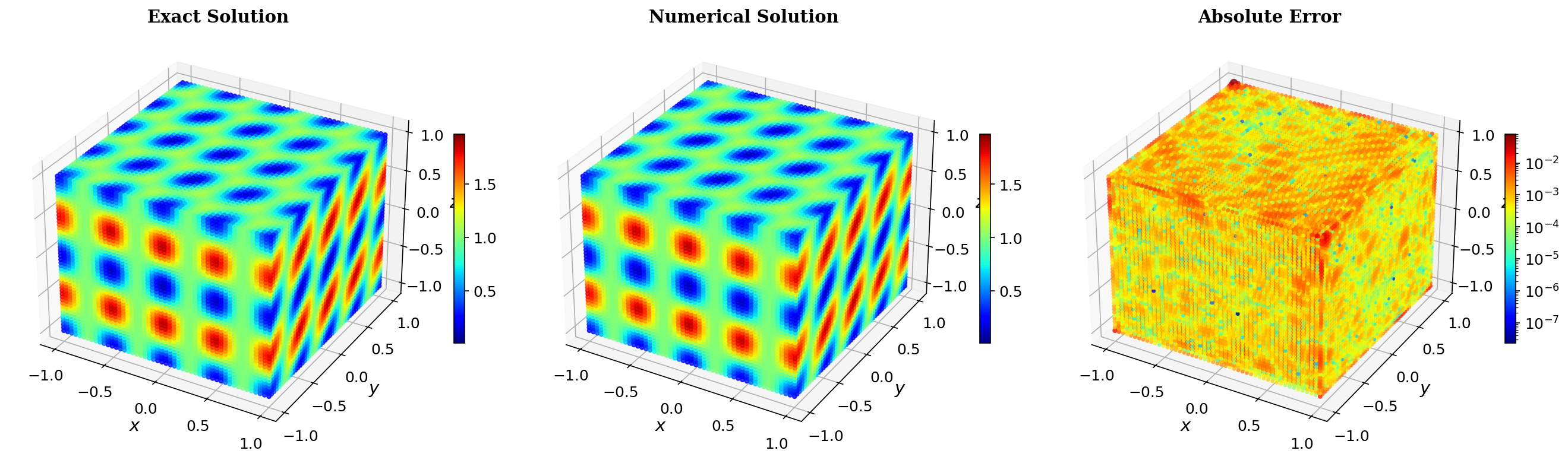}
		\caption{\Cref{ex:3Dfunc2}. Exact function, SCRINN approximation, and pointwise absolute error for the three-dimensional test problem, shown from left to right.}
		\label{fig:3Dfunction}
	\end{figure}

	\Cref{fig:3Dfunction} shows the exact function, the SCRINN approximation, and
	the corresponding pointwise absolute error. The approximation captures the
	main oscillatory and coupled features of the target function, with the
	pointwise error remaining small over most of the computational domain.

	\begin{table}[htbp]
		\centering
		\caption{\Cref{ex:Michalewicz}. Approximation results of different models for the Michalewicz function with $m=10$.}
		\label{tab:3.1.4tab1}
		\begin{tabular}{ccccc}
			\toprule
			Method & Architecture & $E_{L^2}$ & Rank & Training time(s) \\
			\midrule
			
			\multirow{3}{*}{ELM}
			& $[2,512,2000,1]$ &$ 2.9314\times 10^{+1} $& 913 & 2.3 \\
			& $[2,512,4000,1]$ & $1.0077 \times 10^{+2}$ & 1132 & 3.9 \\
			& $[2,512,6000,1]$ & $3.9839 \times 10^{+1}$ & 2414 & 7.7 \\
			
			\midrule
			
			\multirow{3}{*}{RINN-es}
			& $[2,512,2000,1]$ & $ 2.9314\times 10^{+1} $ & 913 & 5.1 \\
			& $[2,512,4000,1]$ & $1.0077 \times 10^{+2}$ & 1132 & 10.3 \\
			& $[2,512,6000,1]$ & $3.9839 \times 10^{+1}$ & 2414 & 19.6 \\
			
			\midrule
			
			\multirow{3}{*}{SCRINN}
			& $[2,2000,2000,1]$ & $1.7179 \times 10^{-3}$ & 890 & 1.5 \\
			& $[2,4000,4000,1]$ & $1.0276 \times 10^{-4}$ & 1621 & 3.6 \\
			& $[2,6000,6000,1]$ & $5.9261 \times 10^{-5}$ & 2933 & 4.4 \\
			%& $[2,8000,8000,1]$ & $7.9193 \times 10^{-5}$ & 3791 & 7.93 \\
			\bottomrule
		\end{tabular}
	\end{table}
	
	\subsection{SCRINN for solving PDEs}
	
	This section evaluates the performance of the proposed SCRINN for solving PDEs. 
	%We first investigate whether progressively increasing the rank of the hidden-layer output matrix during Stage I training can effectively improve the solution accuracy of SCRINN. 
	%We then compare SCRINN with the existing RINN-es model, which has been shown to significantly outperform PIELM \cite{PengHuangYi2026}, as well as its enhanced variant, RINN+, introduced in \Cref{RINN+}.
	
	\subsubsection{SCRINN for steady-state PDEs}
	We begin by considering steady-state PDEs. 
	\begin{example}\label{ex:2DPoisson}
		Consider the two-dimensional Poisson problem
		\begin{equation}\label{poissoneq}
			\begin{cases}
				-\Delta u(x,y)=f(x,y), & (x,y)\in\Omega,\\
				u(x,y)=g(x,y), & (x,y)\in\partial\Omega,
			\end{cases}
		\end{equation}
		where $\Omega=[-1,1]^2$.
		We consider the analytical solution
		\begin{equation} \label{3.3.2liner1}
			u(x,y) = \sin(\pi x) \sin(\pi y) + \sin(4\pi x) \sin(4\pi y) + \sin(8\pi x) \sin(8\pi y).
		\end{equation} 
		This solution contains multiple frequency components and therefore exhibits a multiscale structure, providing a more challenging test for the approximation capability of the proposed method. The source term $f$ and boundary data $g$ are determined accordingly from the analytical solution.
	\end{example}

	\Cref{tab:two-dimensional Poissontab} compares SCRINN and RINN-es in terms of relative $L^2$ error, rank, and training time. SCRINN consistently achieves substantially smaller errors than RINN-es, with the best error being approximately three orders of magnitude lower. Moreover, although the rank generally increases with the number of neural basis functions, a larger rank does not necessarily lead to higher accuracy. For example, SCRINN achieves an error of $8.0784\times10^{-9}$ with rank 1597, whereas RINN-es yields an error of $8.3797\times10^{-3}$ with rank 4010. This indicates that the sparse connectivity and local-support structure of SCRINN, together with the near-orthogonality constraints, lead to a more effective approximation space without requiring a large number of independent basis functions. SCRINN also substantially reduces the computational cost, with training times of only a few seconds for all tested configurations.
	
	\begin{table}[htbp]
		\centering
		\setlength{\tabcolsep}{4pt}
		\renewcommand{\arraystretch}{1.18}
		\caption{\Cref{ex:2DPoisson}. Comparison of approximation accuracy, rank, and computational cost for the two-dimensional Poisson equation.}
		\label{tab:two-dimensional Poissontab}
		\begin{tabular}{ccccc}
			\toprule
			Method & Architecture
			& $E_{L^2}$ & Rank & Time (s)  \\
			\midrule
			
			\multirow{4}{*}{RINN-es}
			& $[2,512,1024,1]$
			& $7.7566\times10^{-1}$ & 731 & 27.2 \\
			
			& $[2,512,2000,1]$
			& $5.2831\times10^{-5}$ & 1021 & 30.7 \\
			
			& $[2,512,4000,1]$
			& $4.1111\times10^{-6}$ & 1971 & 43.25 \\
			
			& $[2,512,6000,1]$
			& $8.3797\times10^{-3}$ & 4010 & 26.9 \\
			
			\cmidrule(lr){1-5}
			
			\multirow{3}{*}{SCRINN}
			& $[2,2000,2000,1]$
			& $1.2807\times10^{-5}$ & 966 & 3.51 \\
			
			& $[2,4000,4000,1]$
			& $4.8480\times10^{-8}$ & 1516 & 4.51 \\
			
			& $[2,6000,6000,1]$
			& $8.0784\times10^{-9}$ & 1597 & 6.7 \\
			
			\bottomrule
		\end{tabular}
	\end{table} 
	
	\begin{table}[htpb]
		\centering
		\caption{\Cref{ex:2DPoisson}. Influence of different loss functions and $\varepsilon$ values on the approximation accuracy of RINN-es and SCRINN. Boldface entries correspond to the results obtained by RINN+.}
		\label{tab:lossdiif}
		\begin{tabular}{ccccc}
			\toprule
			\text{Method} & \text{Architecture} &
			$\boldsymbol{\varepsilon}$ &
			\text{Original Loss} \eqref{lossf} &
			\text{Proposed Loss} \eqref{lossscrinn} \\
			\midrule
			
			\multirow{3}{*}{RINN-es}
			& \multirow{3}{*}{$[2,512,4000,1]$}
			& 0.01 & $4.1110\times10^{-6}$ & $\mathbf{3.9364\times10^{-6}}$ \\
			&& 0.5 & $1.9069\times10^{-3}$ & $\mathbf{4.1882\times10^{-6}}$ \\
			&& 1.0 & $3.6715\times10^{-3}$ & $\mathbf{4.1565\times10^{-6}}$ \\
			
			\cmidrule{1-5}
			
			\multirow{3}{*}{SCRINN}
			& \multirow{3}{*}{$[2,4000,4000,1]$}
			& 0.01 & $1.0430\times10^{-7}$ & $1.5351\times10^{-7}$ \\
			&& 0.5 & $1.0731\times10^{-8}$ & $6.7763\times10^{-8}$ \\
			&& 1.0 & $1.3524\times10^{-7}$ & $4.8480\times10^{-8}$ \\
			
			\bottomrule
		\end{tabular}
	\end{table}

	\begin{figure}[htbp]
		\centering
		\includegraphics[width=0.85\linewidth]{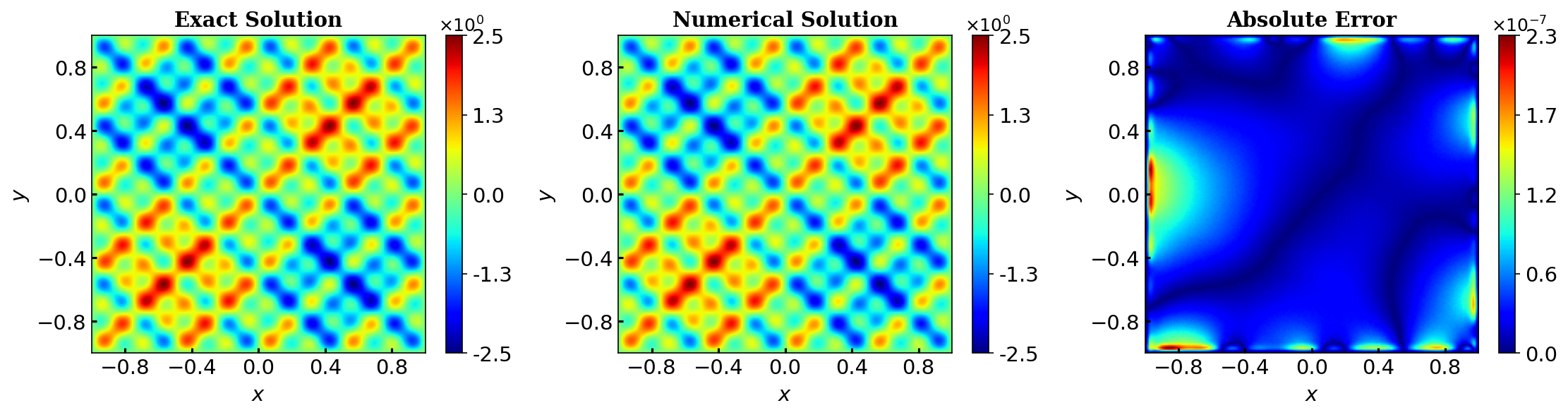}
		\caption{\Cref{ex:2DPoisson}. Results of the model for the two-dimensional Poisson equation. Each column shows: (left) the exact solution, (middle) the predicted solution, and (right) the point-wise absolute error. }
		\label{fig:posong_comparisonz}
	\end{figure}
	
	To examine the sensitivity to the hyperparameter $\varepsilon$, we compare the original and proposed loss functions in \Cref{tab:lossdiif}. For RINN-es, the original loss exhibits strong sensitivity to $\varepsilon$: the relative $L^2$ error increases from $4.11\times10^{-6}$ at $\varepsilon=0.01$ to $3.67\times10^{-3}$ at $\varepsilon=1$. In contrast, the proposed loss maintains errors of approximately $4\times10^{-6}$ across all tested values. Thus, the proposed loss substantially improves the robustness of RINN-es with respect to $\varepsilon$. For SCRINN, both losses remain robust, with errors consistently within the range of $10^{-8}$--$10^{-7}$. This further indicates that SCRINN is relatively insensitive to $\varepsilon$, while the proposed loss preserves its high approximation accuracy.

\begin{figure}[htbp]
		\centering
		\begin{subfigure}{0.45\textwidth}
			\centering
			\includegraphics[width=0.9\linewidth]{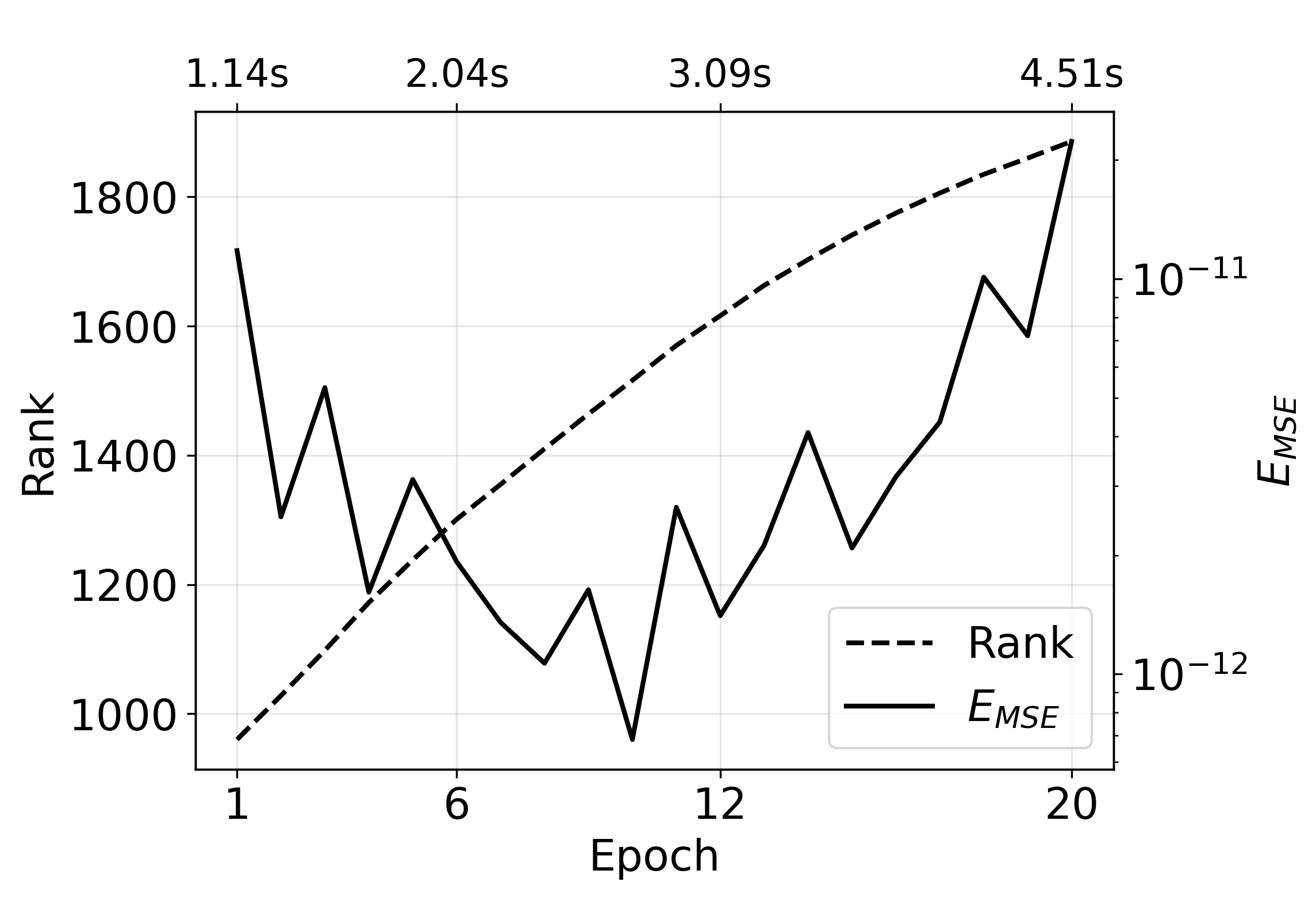}
		\end{subfigure}
		\hfill
		\begin{subfigure}{0.45\textwidth}
			\centering
			\includegraphics[width=0.9\linewidth]{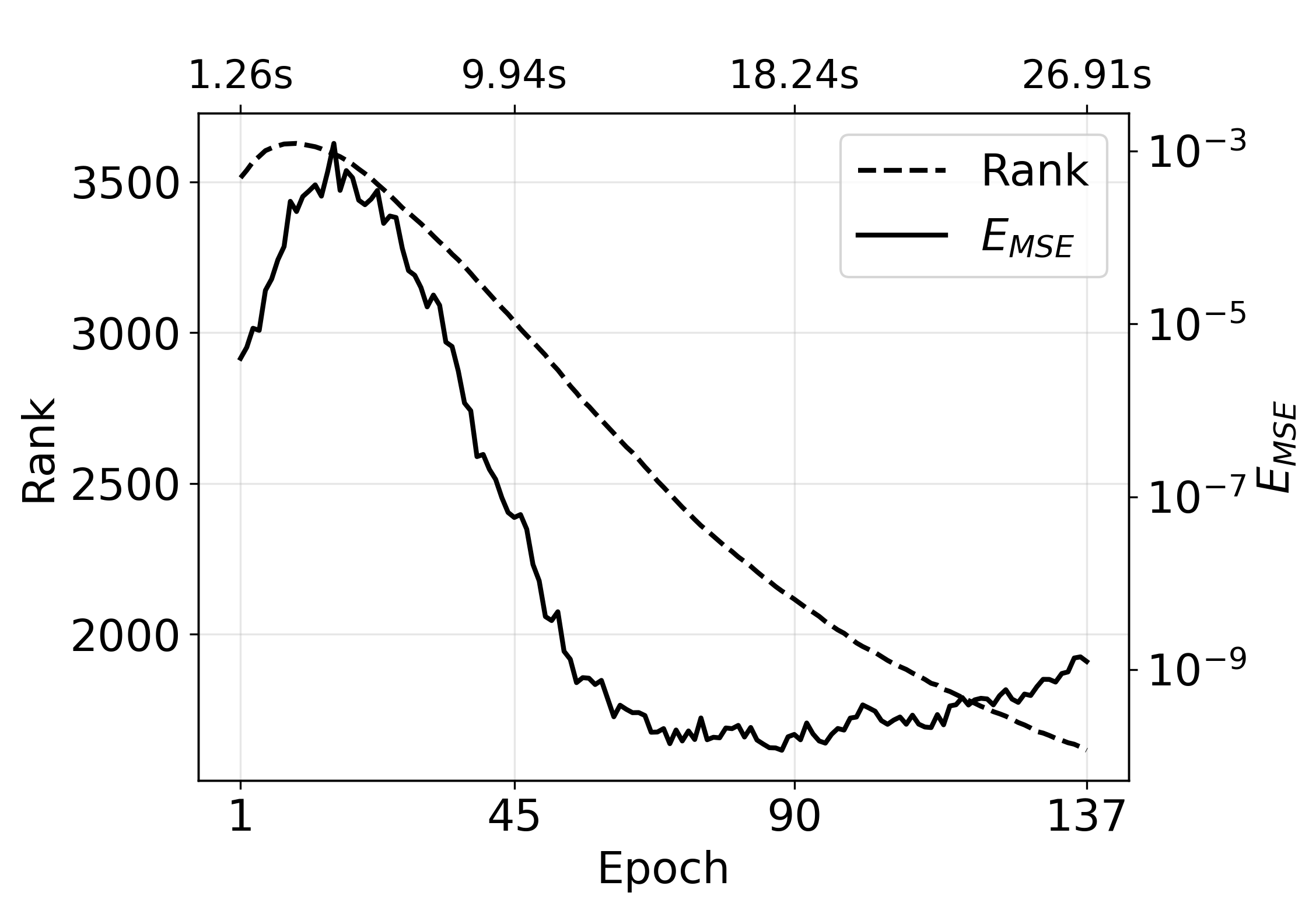}
		\end{subfigure}  
		\caption{\Cref{ex:2DPoisson}. $E_{\mathrm{MSE}}$ and rank evolution with $4000$ neural basis functions. Left: SCRINN; Right: RINN-es.}
		\label{fig:posong3.3.2}
	\end{figure}

	Finally, we examine the evolution of the MSE-based error $E_{\mathrm{MSE}}$ and the rank of the neural basis function matrix during Stage I training. As shown in \Cref{fig:posong3.3.2}, the rank of SCRINN generally increases as $E_{\mathrm{MSE}}$ decreases and remains at a low level, indicating progressive enrichment of the neural basis space. In contrast, RINN-es exhibits an overall decrease in rank during training, accompanied by larger variations in $E_{\mathrm{MSE}}$. These results suggest that the sparse connectivity structure of SCRINN promotes a more effective utilization of the learned basis functions and provides a more favorable relationship between rank enrichment and error reduction.
	
	% Finally, we investigate the relationship between the MSE-based error $E_{\mathrm{MSE}}$ and the rank of the neural basis function matrix in SCRINN and RINN-es. \Cref{fig:posong3.3.2} illustrates the evolution of these two quantities during the Stage I training process.
	% The training dynamics further demonstrate the advantages of the SCRINN architecture. As shown in \Cref{fig:posong3.3.2}, during SCRINN training, the rank of the neural basis function matrix continuously increases, while the error $E_{\mathrm{MSE}}$ decreases correspondingly and remains at a very low level. This behavior indicates that the sparse connectivity pattern effectively promotes the enrichment of the neural basis space, providing additional independent degrees of freedom for reducing the PDE residual. In contrast, for RINN-es, the rank of the neural basis function matrix exhibits an overall decreasing trend during training, while the reduction of $E_{\mathrm{MSE}}$ is accompanied by noticeable changes. This suggests that the fully connected architecture may suffer from ineffective utilization of the available basis functions, where the increase in network parameters does not necessarily transfer to enhanced approximation capability. These observations demonstrate that SCRINN establishes a more favorable coupling between rank evolution and loss reduction through its structural design, which is one of the key factors enabling it to achieve highly accurate solutions.
	
	\subsubsection{SCRINN for time-dependent PDEs}
	Next, we investigate the effectiveness and robustness of SCRINN for solving time-dependent PDEs, including the
	advection equation, heat equation, and wave equation.
	
	\begin{example}\label{ex:1advection}
		
		We first consider the one-dimensional advection equation on the spatial domain
		$\Omega=[-1,1]$ and the time interval $I=[0,1]$:
		\begin{equation}\label{eq:advection}
			\begin{cases}
				u_t(x,t) = -c u_x(x,t), & (x,t)\in\Omega\times I,\\
				u(-1,t)=u(1,t), & t\in I,\\
				u(x,0)=h(x), & x\in\Omega.
			\end{cases}
		\end{equation}
		Here, $c=0.4$ is the wave speed, and periodic boundary conditions are imposed.
		The initial condition is given by
		\begin{equation}\label{3.6.1liner1}
			h(x)=e^{0.5\sin(2\pi x)}-1.
		\end{equation}
	\end{example}
	
\begin{figure}[htbp]
		\centering
		\includegraphics[width=0.85\linewidth]{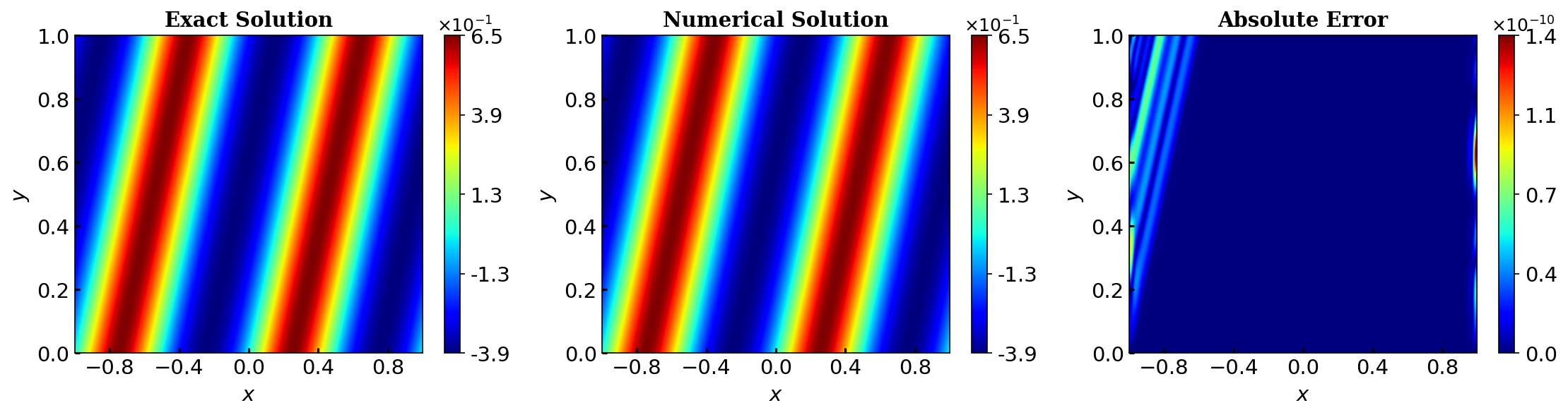}   
		\caption{\Cref{ex:1advection}. Results of the model for the one-dimensional advection equation. From left to right, the figure shows the exact function, the SCRINN approximation, and pointwise absolute error.}
		\label{fig:advectionjie}
	\end{figure}

	\Cref{tab:advectiontab} compares the performance of RINN-es and SCRINN for solving the advection equation \eqref{eq:advection} with the initial condition \eqref{3.6.1liner1}. As shown in \Cref{tab:advectiontab}, SCRINN consistently outperforms RINN-es in terms of both approximation accuracy and computational efficiency. Specifically, SCRINN reduces the relative $L^2$ error from $2.1387\times10^{-9}$ to $3.0506\times10^{-11}$ while decreasing the training time from $14.23$ s to $5.41$ s. Moreover, SCRINN achieves a substantially higher effective rank ($1200$ versus $717$), indicating a stronger expressive capability of the learned representation. The corresponding predicted solution and the pointwise absolute error are shown in \Cref{fig:advectionjie}, further demonstrating the excellent approximation performance of SCRINN.

	\begin{table}[htbp]
		\centering
		\caption{\Cref{ex:1advection}. Comparison of RINN-es and SCRINN for the one-dimensional advection equation with the initial condition \eqref{3.6.1liner1}.}
		\label{tab:advectiontab}
		
		\begin{tabular}{ccccc}
			\toprule
			Method & Architecture & $E_{L^2}$ & Rank & Time (s) \\
			\midrule
			RINN-es & $[2,512,4000,1]$ & $2.1387 \times 10^{-9}$ & 717 & 14.23 \\
			SCRINN  & $[2,4000,4000,1]$ & $3.0506 \times 10^{-11}$ & 1200 & 5.41 \\
			\bottomrule
		\end{tabular}
		
	\end{table}

	\begin{example}\label{ex:1aheat}
		Next, we consider the heat conduction equation on $\Omega=[-1,1]$ over
		$I=[0,1]$:
		\begin{equation}\label{eq:heat}
			\begin{cases}
				u_t(x,t)-u_{xx}(x,t)=f(x,t), &(x,t)\in\Omega\times I,\\
				u(x,t)=g(t), &(x,t)\in\partial\Omega\times I,\\
				u(x,0)=h(x), &x\in\Omega.
			\end{cases}    
		\end{equation}
		The exact solution is prescribed as
		\begin{equation}\label{3.7.1liner1}
			u(x,t)=e^{-t}\sin(k\pi x),\qquad k=6,12,
		\end{equation}
		and the source term and boundary data are obtained accordingly. Here, the parameter $k$ is the spatial frequency mode number (also called the wave number or Fourier mode index).
	\end{example}   
	
\begin{figure}[htbp]
		\centering
		\includegraphics[width=0.85\linewidth]{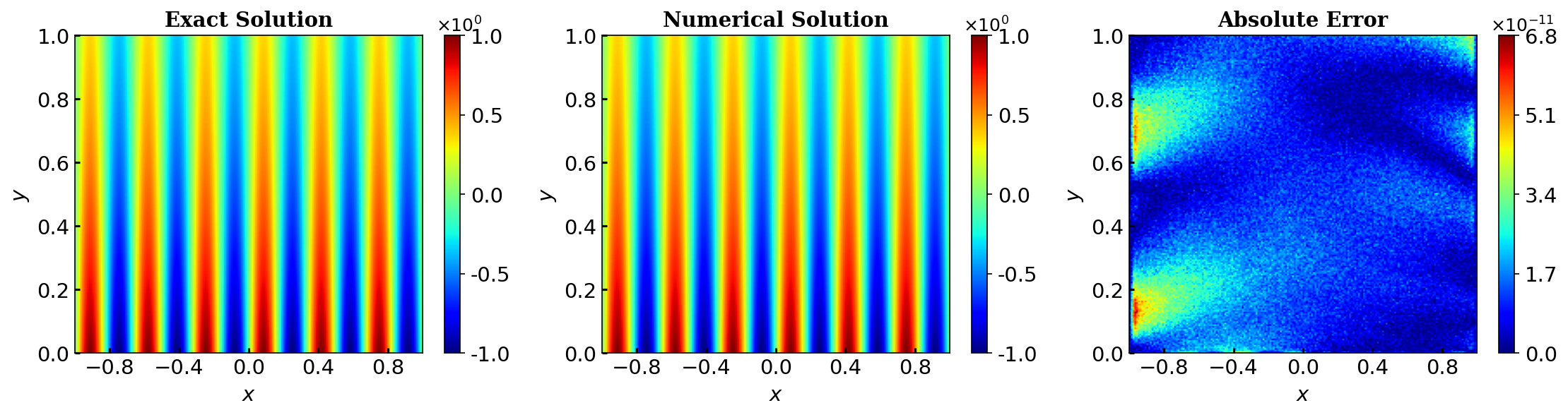} 
		\includegraphics[width=0.85\linewidth]{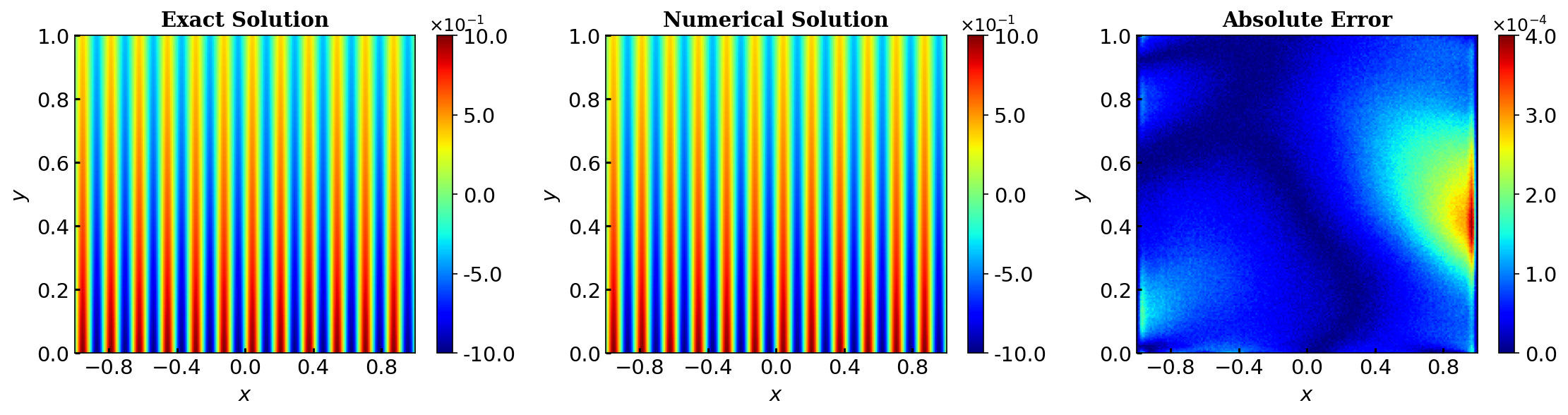} 
		\caption{\Cref{ex:1aheat}. Approximation of SCRINN for heat equation \eqref{eq:heat}. From left to right, the columns show the exact function, the SCRINN approximation, and pointwise absolute error. The rows from top to bottom correspond to the results for $k=6$ and $k=12$ in \eqref{3.7.1liner1} , respectively.}
		\label{fig:heatjie}
	\end{figure}

	\Cref{tab:heattab} compares the performance of RINN-es and SCRINN for the one-dimensional heat conduction equation under different values of $k$. The table reports the relative $L^2$ error, the rank of the neural basis functions, and the corresponding training time under the optimized parameter settings. As $k$ increases, both methods exhibit a deterioration in accuracy; however, SCRINN consistently achieves lower relative $L^2$ errors than RINN-es. In particular, for $k=6$, SCRINN reduces the relative $L^2$ error from $3.65\times10^{-7}$ to $3.25\times10^{-11}$ while requiring only $7.72$ seconds of training time, compared with $13.00$ seconds for RINN-es. For $k=12$, SCRINN likewise achieves a lower error of $1.86\times10^{-4}$ compared with $2.06\times10^{-3}$ for RINN-es, while reducing the training time from $13.44$ to $4.78$ seconds. These results demonstrate that SCRINN provides a favorable balance between approximation accuracy and computational efficiency. \Cref{fig:heatjie} presents the corresponding numerical solutions and error distributions.
	
	\begin{table}[htbp]
		\centering
		\caption{\Cref{ex:1aheat}. Comparison results for the one-dimensional heat conduction equation}
		\label{tab:heattab}
		\begin{tabular}{cccccc}
			\toprule
			$k$ value & Method & Architecture & $E_{L^2}$ & Rank & Training time(s) \\
			\midrule
			\multirow{2}{*}{6} & RINN-es & $[2,512,4000,1]$ & $3.6575 \times 10^{-7}$ & 573  & 13.00\\
			& SCRINN & $[2,4000,4000,1]$ & $3.2536 \times 10^{-11}$ & 796 & 7.72\\
			\cmidrule{1-6}
			\multirow{2}{*}{12} & RINN-es & $[2,512,4000,1]$ & $2.0613 \times 10^{-3}$ & 645 & 13.44\\
			& SCRINN & $[2,4000,4000,1]$ & $1.8650 \times 10^{-4}$ & 724 & 4.78\\
			\bottomrule
		\end{tabular}
	\end{table}

	\begin{example}\label{ex:1awave}
		Finally, we consider the wave equation on $\Omega=[0,1]$ and $I=[0,1]$:
		\[
		\begin{cases}
			u_{tt}(x,t)=c^2u_{xx}(x,t), &(x,t)\in\Omega\times I,\\
			u(0,t)=u(1,t)=0, &t\in I,\\
			u(x,0)=h(x), &x\in\Omega,\\
			u_t(x,0)=0, &x\in\Omega.
		\end{cases}
		\]
		The wave speed is set as $c=2$. The initial displacement is chosen as
		\begin{equation}\label{3.8.1liner2}
			h(x)=\sin(2\pi x)+\sin(4\pi x).
		\end{equation}
		The corresponding exact solution is
		\[
		u(x,t)
		=
		\sin(2\pi x)\cos(4\pi t)
		+\sin(4\pi x)\cos(8\pi t).
		\]
		
	\end{example}
	
\begin{figure}[htbp]
		\centering
		\includegraphics[width=0.85\linewidth]{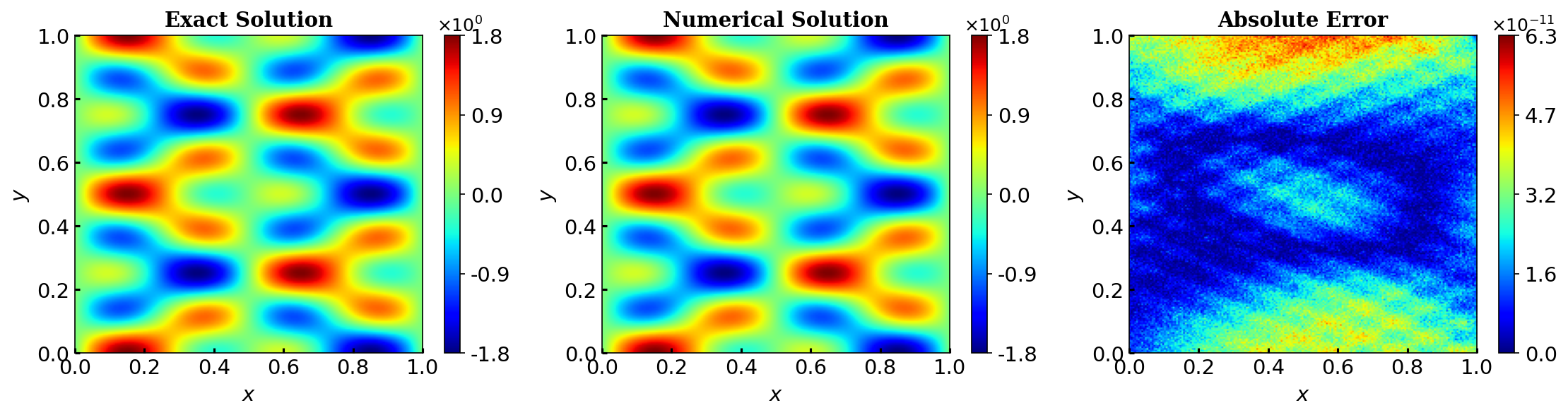} 
		\caption{\Cref{ex:1awave}. Results of the model for the one-dimensional wave equation. From left to right, the figure shows the exact function, the SCRINN approximation, and pointwise absolute error.}
		\label{fig:bodojie}
	\end{figure}

	\Cref{tab:bodotab} summarizes the numerical results for the one-dimensional wave equation under the optimal parameter settings. The table reports the relative $L^2$ errors between the predicted and exact solutions, the ranks of the resulting neural basis functions, and the corresponding training time. It can be observed that SCRINN achieves substantially smaller relative $L^2$ errors and shorter training time than RINN-es. In particular, SCRINN reduces the training time by more than a factor of three while achieving error as low as $3.1039\times 10^{-11}$. \Cref{fig:bodojie} presents the corresponding visualization results.
	
	\begin{table}[htbp]
		\centering
		\caption{\Cref{ex:1awave}. Comparison results for the one-dimensional wave equation.}
		\label{tab:bodotab}
		\begin{tabular}{ccccc}
			\toprule
			Method & Architecture & $E_{L^2}$ & Rank & Training time (s) \\
			\midrule
			RINN-es & $[2,512,4000,1]$ & $7.8580 \times 10^{-8}$ & 655 & 15.30 \\
			SCRINN & $[2,4000,4000,1]$ & $3.1039 \times 10^{-11}$ & 771 & 5.77 \\
			\bottomrule
		\end{tabular}
	\end{table}
	
	\begin{example}\label{ex:3DPoisson}
		We consider the three-dimensional Poisson problem (4.4) on
		$\Omega=[-1,1]^3$ with the analytical solution
		\begin{equation}\label{3.5.2liner1}
			u(\boldsymbol{x})
			=-e^{-x}\sin(y)\cosh(z).
		\end{equation}
		The source term $f$ and boundary data $g$ are determined from the
		analytical solution.
	\end{example}
	
	\Cref{tab:three-ceshitab} compares RINN-es and SCRINN under different
	network architectures and collocation settings. In addition to the default
	configuration, we consider two refined collocation sets with $5832$ and
	$8000$ interior points and $7998$ and $9996$ boundary points, respectively.
	For both methods, increasing the network width can substantially reduce the
	relative $L^2$ error. For example, with $5832$ interior points, increasing
	the number of basis functions from $4000$ to $6000$ reduces the error of
	SCRINN from $8.8023\times10^{-8}$ to $8.2511\times10^{-9}$. With $6000$
	basis functions, further increasing the number of interior points from
	$5832$ to $8000$ reduces the error to $1.5351\times10^{-9}$. The effect of
	collocation refinement is less consistent for RINN-es, for which the error
	increases from $1.4759\times10^{-8}$ to $2.8185\times10^{-8}$ under the two
	finest configurations.
	
	SCRINN also requires substantially less training time than RINN-es across
	all tested configurations. Its training time ranges from $6.23$ to
	$20.10$ seconds, compared with $608.51$ to $801.11$ seconds for RINN-es.
	For the finest configuration, SCRINN achieves a relative $L^2$ error of
	$1.5351\times10^{-9}$ in $20.10$ seconds, whereas RINN-es obtains
	$2.8185\times10^{-8}$ in $801.11$ seconds. Thus, for this problem, the
	sparse architecture provides both lower approximation errors and
	substantially lower training costs.
	
	\begin{table}[htbp]
		\centering
		\renewcommand{\arraystretch}{1.15}
		\caption{\Cref{ex:3DPoisson}. Comparison results for the three-dimensional Poisson equation.}
		\label{tab:three-ceshitab}
		\begin{tabular}{cclcccc}
			\toprule
			Method & Architecture & $K_{\rm res}$ & $K_{\rm bcs}$
			& $E_{L^2}$ & Rank & Time(s) \\
			\midrule
			
			\multirow{4}{*}{RINN-es}
			& $[2,512,4000,1]$
			& 4096 & 6441
			& $2.1679 \times 10^{-7}$ & 3880 & 608.51 \\
			
			& $[2,512,4000,1]$
			& 5832 & 7998
			& $4.4513 \times 10^{-7}$ & 4000 & 678.01 \\
			
			& $[2,512,6000,1]$
			& 5832 & 7998
			& $ 1.4759 \times 10^{-8}$ & 4096 & 718.62 \\
			
			& $[2,512,6000,1]$
			& 8000 & 9996
			& $2.8185 \times 10^{-8}$ & 5895 & 801.11 \\
			
			\cmidrule{1-7}
			
			\multirow{4}{*}{SCRINN}
			& $[2,4000,4000,1]$
			& 4096 & 6441
			& $1.0067 \times 10^{-7}$ & 3745 & 6.23 \\
			
			& $[2,4000,4000,1]$
			& 5832 & 7998
			& $8.8023 \times 10^{-8}$ & 3999 & 11.30 \\
			
			& $[2,6000,6000,1]$
			& 5832 & 7998
			& $8.2511 \times 10^{-9}$ & 3881 & 15.80 \\
			
			& $[2,6000,6000,1]$
			& 8000 & 9996
			& ${1.5351 \times 10^{-9}}$ & 4911 & 20.10 \\
			
			\bottomrule
		\end{tabular}
	\end{table}

\begin{figure}[htbp]
		\centering
		\includegraphics[width=0.85\linewidth]{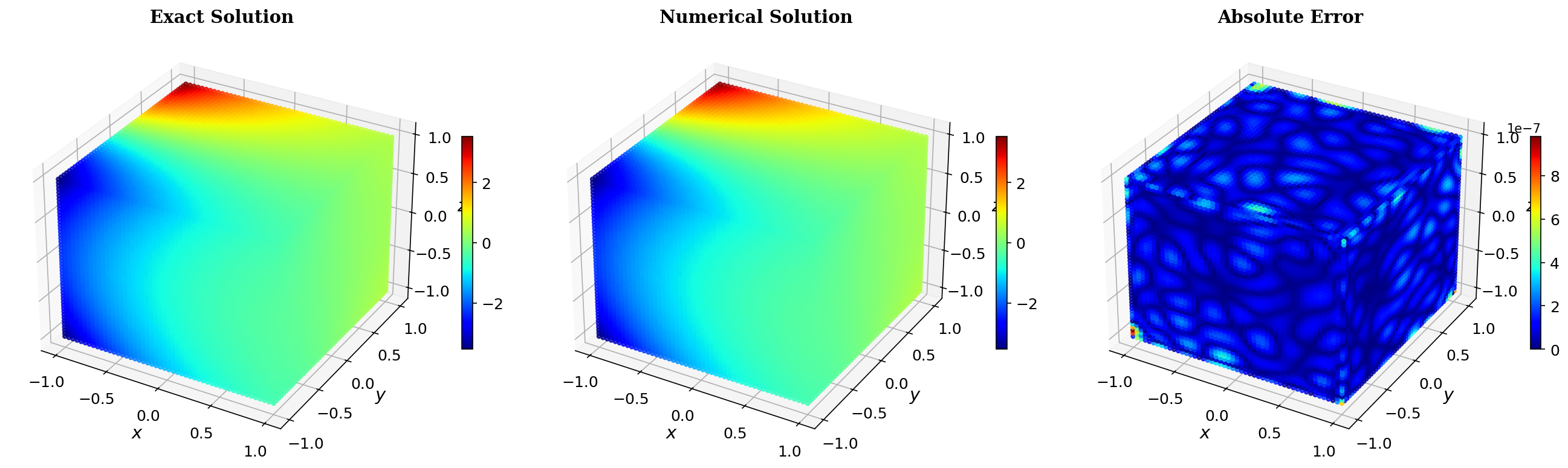}
		\caption{\Cref{ex:3DPoisson}. Visualization of the numerical solution for the three-dimensional Poisson equation. From left to right, the figures display the analytical solution, the SCRINN prediction, and the pointwise absolute error.}
		\label{fig:3djie}
	\end{figure}

	\Cref{fig:3djie} shows the analytical solution, the SCRINN approximation,
	and the corresponding pointwise absolute error. The numerical solution
	captures the main features of the analytical solution, with small pointwise
	errors over the computational domain.

\begin{figure}[htbp]
		\centering
		\begin{subfigure}{0.45\textwidth}
			\centering
			\includegraphics[width=0.9\linewidth]{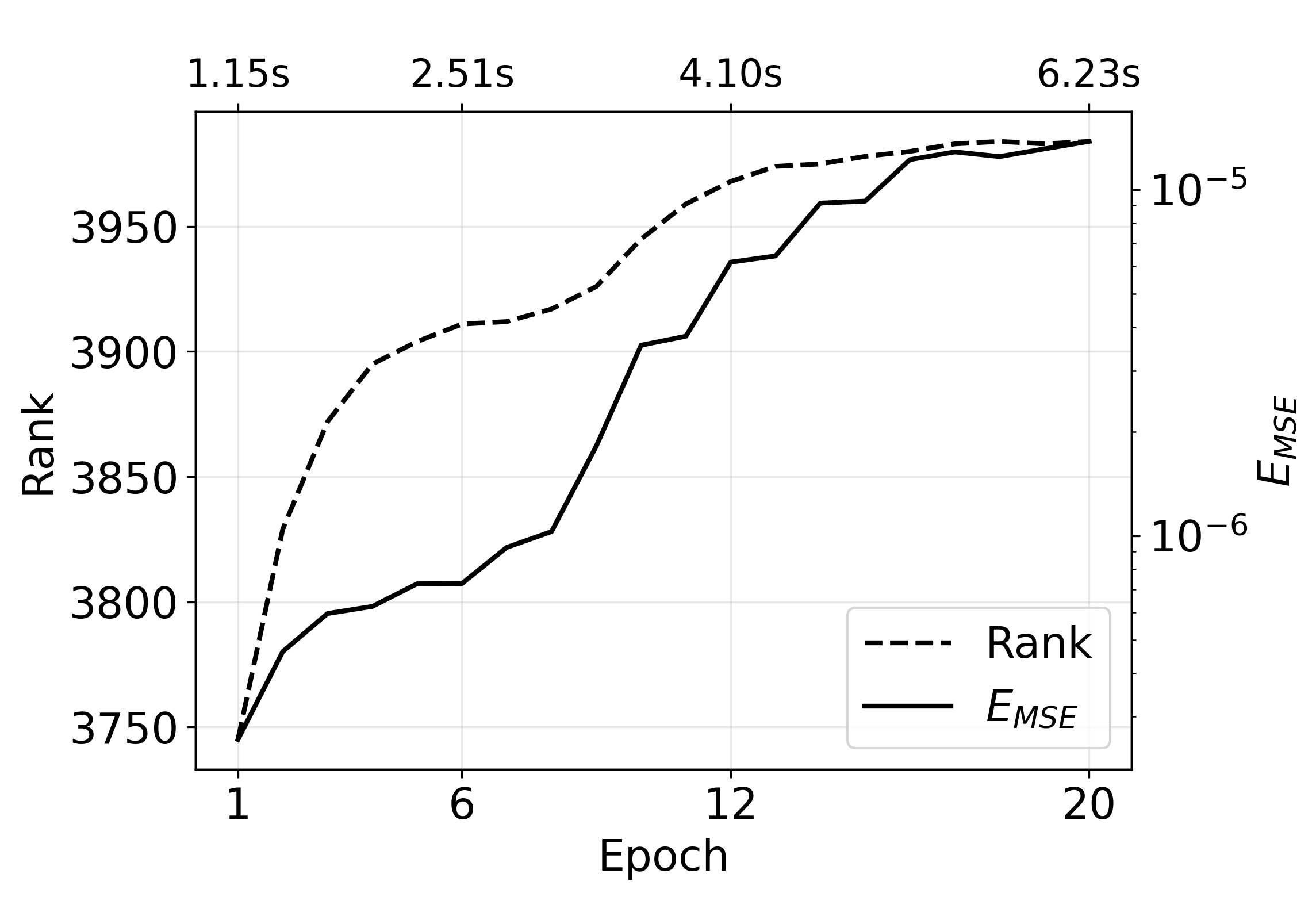}
		\end{subfigure}
		\hfill
		\begin{subfigure}{0.45\textwidth}
			\centering
			\includegraphics[width=0.9\linewidth]{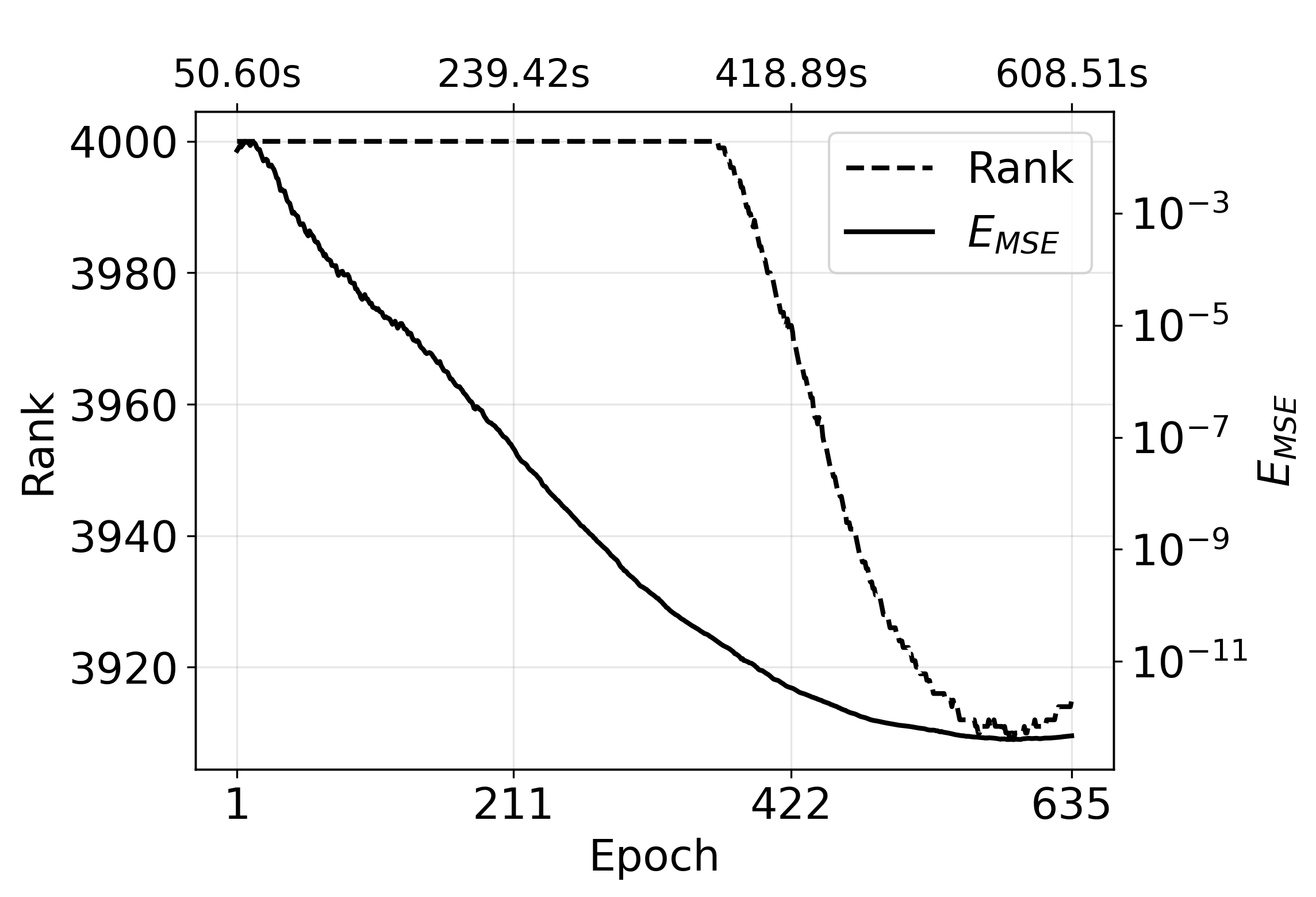}
		\end{subfigure}
		\caption{\Cref{ex:3DPoisson}. $E_{\mathrm{MSE}}$ and rank evolution for the  three-dimensional Poisson equation with $4000$ neural basis functions. Left: SCRINN; Right: RINN-es.}
		\label{fig:posong_comparison}
	\end{figure}

	Finally, we examine the evolution of the MSE-based error
	$E_{\mathrm{MSE}}$ and the rank of the neural basis function matrix during
	Stage I training with $4000$ basis functions; see
	\Cref{fig:posong_comparison}. The effective ranks of both SCRINN and RINN-es
	remain close to their respective full ranks throughout the training process.
	This indicates that both methods maintain highly expressive neural bases. Although RINN-es continues to reduce
	the training residual with prolonged optimization, the resulting relative
	$L^2$ error remains larger than that of SCRINN, as reported in
	\Cref{tab:three-ceshitab}.
	
	\subsubsection{SCRINN for nonlinear PDEs}

	\begin{example}\label{ex:2ha}
		We consider a nonlinear equation with Dirichlet boundary conditions in $\Omega=[-1,1]^2$,
		\begin{equation}\label{NonLEQ}
			\left\{\begin{aligned} 
				\Delta u - 5u + \sin(u) = f, &\qquad  x \in \Omega, \\
				u = g, &\qquad  x \in \partial \Omega.
			\end{aligned}\right.        
		\end{equation}
		% The equation is defined on the square domain $\Omega = [-1, 1]^2$ and has 
		Consider the following analytical solution:
		\begin{equation} \label{3.4.2liner1}
			\begin{split}
				u(\boldsymbol{x}) =& \sin(2\pi x) \cos(\pi x) + \cos(3\pi y) \sin(1.5\pi y) \\
				&+ 0.5x^3 - 0.3x^2 + 0.2y^3 - 0.1y^2 + 0.8xy + 0.4x^2y - 0.2xy^2 \\
				&+ 0.3 \tanh(2x) \sinh(y) + 0.2 e^{-(x^2 + y^2)/2} + 1.0
			\end{split}
		\end{equation}
		
	\end{example}
	
\begin{figure}[htbp]
		\centering
		\begin{subfigure}{0.45\textwidth}
			\centering
			\includegraphics[width=0.9\linewidth]{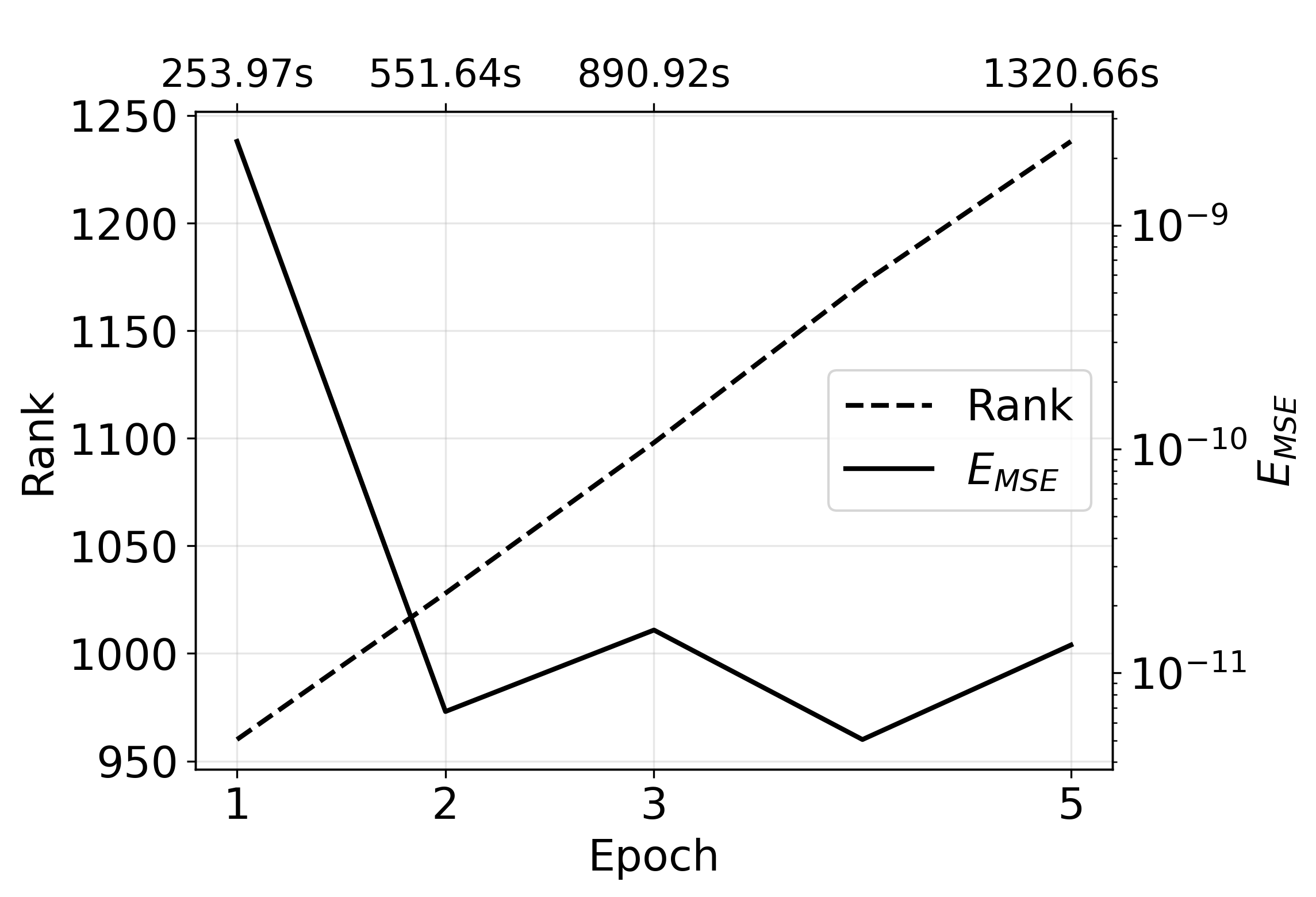} 
		\end{subfigure}
		\hfill
		\begin{subfigure}{0.45\textwidth}
			\centering
			\includegraphics[width=0.9\linewidth]{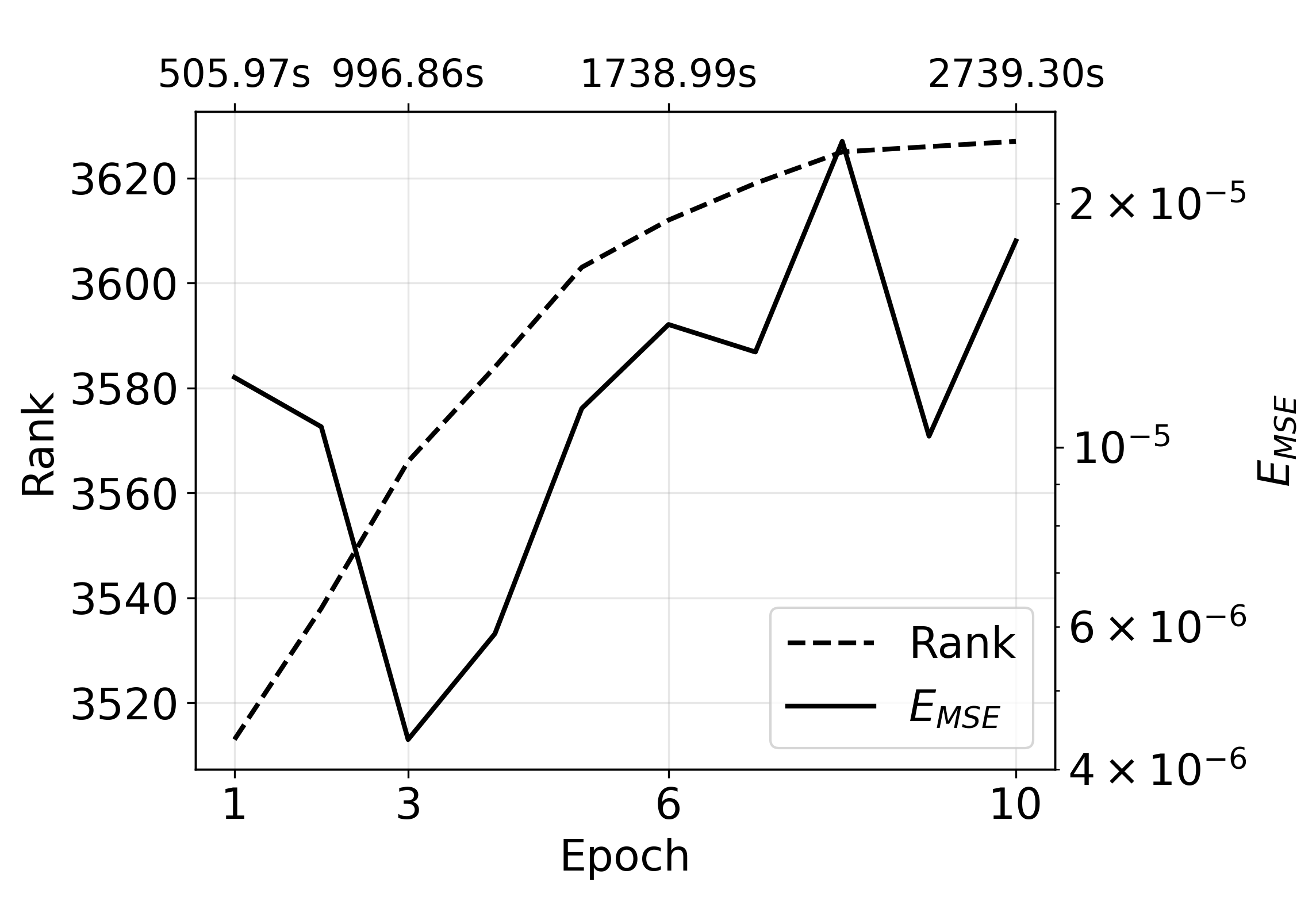} 
		\end{subfigure}  
		\caption{\Cref{ex:2ha}. $E_{\mathrm{MSE}}$ and rank evolution of SCRINN and RINN-es for the two-dimensional nonlinear Helmholtz equation with $4000$ neural basis functions. Left: SCRINN; Right: RINN-es.}
		\label{fig:heat3}
	\end{figure}

	First, we examine the evolution of the MSE-based error \(E_\mathrm{MSE}\) and the rank of the neural basis function matrix. For SCRINN, we set the maximum number of epochs \(\nu_{\max}=5\) in Stage I, while for RINN-es, it is set to \(10\). 
	\Cref{fig:heat3} presents the evolution curves of these two metrics during model training.
	As shown in \Cref{fig:heat3}, the steady growth of the rank of the neural basis functions exhibits strong synchronization with the continuous decrease of the MSE-based error $E_{\mathrm{MSE}}$. The training dynamics indicate that when handling the nonlinear problem, the model reduces the error to a low level after the first optimization stage. This confirms that the sparse connection architecture effectively reduces the model’s dependence on parameter initialization. The underlying mechanism lies in the fact that the sparse structure imposes a prior constraint of local support on each neuron, naturally restricting the outputs of randomly initialized neurons to specific subregions of the input space. Benefiting from this structural constraint, the model avoids the high correlation among basis functions caused by global random connections in fully connected networks, thereby enabling the formation of a well-conditioned initial representation space with low redundancy even at the early stage of training.
	
\begin{figure}[htbp]
		\centering
		\includegraphics[width=0.85\linewidth]{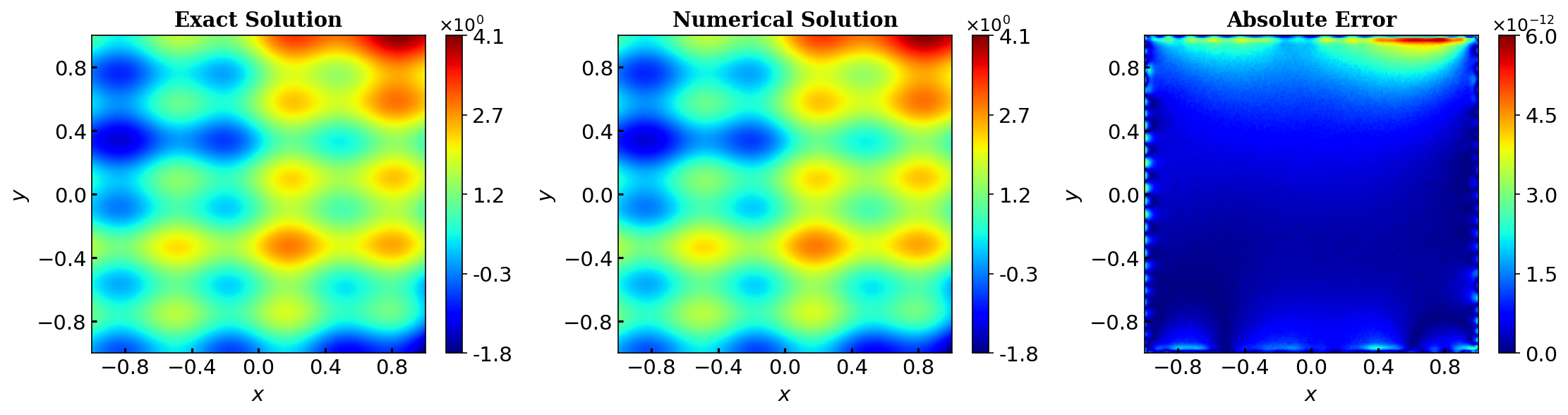}
		\caption{\Cref{ex:2ha}. Approximation results of SCRINN for the nonlinear equation \eqref{NonLEQ}. 
			From left to right, the figure shows the exact function, the SCRINN approximation, and
			pointwise absolute error.}
		\label{fig:haimujie}
	\end{figure}

	\Cref{tab:Helmholtztab} summarizes the numerical results for the equation \eqref{NonLEQ}. For each analytical solution, we report the relative $L^2$ error, the rank of the neural basis functions, and the training time obtained under the optimal parameter settings. The results demonstrate that SCRINN consistently achieves substantially higher accuracy with a much lower rank and shorter training time than RINN-es. In particular, SCRINN reduces the rank by more than 70\% while simultaneously reducing the relative $L^2$ error by several orders of magnitude. These results indicate that the sparse connectivity of SCRINN enables a more efficient representation of the solution space without requiring a large number of neural basis functions. The corresponding numerical solutions and error distributions are visualized in \Cref{fig:haimujie}.

	% \begin{table}[htbp]
		% \centering
		% \caption{\Cref{ex:2ha}. Comparison results for the equation \eqref{NonLEQ}}
		% \label{tab:Helmholtztab}
		% \begin{tabular}{cccccc}
			% \toprule
			% Analytical solution & Method & Architecture & $E_{L^2}$ & Rank & Training time(s) \\
			% \midrule
			% \multirow{2}{*}{\eqref{3.4.1liner1}} & RINN-es & $[2,512,4000,1]$ & $1.5281 \times 10^{-6}$ & 3538 & 5374.07 \\
			% & SCRINN & $[2,4000,4000,1]$ & $2.2711 \times 10^{-12}$ & 960 & 3459.45 \\
			% \cmidrule{1-6}
			% \multirow{2}{*}{\eqref{3.4.2liner1}} & RINN-es & $[2,512,4000,1]$ & $1.9482 \times 10^{-8}$ & 3566 & 2739.30 \\
			% & SCRINN & $[2,4000,4000,1]$ & $7.6778 \times 10^{-13}$ & 1172 & 1320.66 \\
			% \bottomrule
			% \end{tabular}
		% \end{table}
	
	\begin{table}[htbp]
		\centering
		\caption{\Cref{ex:2ha}. Comparison results for the equation \eqref{NonLEQ}}
		\label{tab:Helmholtztab}
		\begin{tabular}{ccccc}
			\toprule
			Method & Architecture & $E_{L^2}$ & Rank & Time (s) \\
			\midrule
			RINN-es & $[2,512,4000,1]$ & $1.9482 \times 10^{-8}$ & 3566 & 2739.30 \\
			SCRINN & $[2,4000,4000,1]$ & $7.6778 \times 10^{-13}$ & 1172 & 1320.66 \\
			\bottomrule
		\end{tabular}
	\end{table}

	\section{Conclusion}\label{conclusion}
	
	In this work, we proposed SCRINN, a sparsely connected rank-inspired neural network for function approximation and the numerical solution of PDEs. By combining channel-wise sparse connectivity with a Gram-matrix-based orthogonality- and normalization-guided loss, SCRINN constructs a compact neural basis with sustained rank growth while reducing the computational cost of full connectivity. The proposed design also alleviates the sensitivity to the hyperparameter $\varepsilon$ observed in RINN loss function.
	
	Numerical experiments on function approximation and a range of steady-state, time-dependent, and nonlinear PDEs demonstrate that SCRINN achieves comparable or higher accuracy than RINN-es with substantially lower training cost, often using a lower-rank neural basis. These results highlight the importance of basis quality and effective utilization rather than simply increasing the number of basis functions. We also showed that the proposed loss improves the robustness of the original RINN architecture, leading to the RINN+ variant. Future work will focus on more flexible sparse connectivity and problem-adaptive basis construction for higher-dimensional and more challenging PDEs.

	\bibliographystyle{siamplain}
	\bibliography{references}
\end{document}

%% file: ex_shared.tex
\usepackage{lipsum}
\usepackage{amsfonts}
\usepackage{graphicx}
\usepackage{epstopdf}
\usepackage{algorithm}
\usepackage{algorithmic}
\usepackage{mathrsfs}

\usepackage{amsmath}\allowdisplaybreaks  % 允许跨页显示公式
\usepackage{amssymb}        % 数学符号
\usepackage{mathrsfs}       % 特殊花体符号
\usepackage{cases}          % 分段函数
\usepackage{siunitx}        % 物理单位和科学计数
\usepackage{amsmath} 
\usepackage{array}
\usepackage{booktabs}       % 专业表格
\usepackage{longtable}      % 跨页表格
\usepackage{multirow}
\usepackage{makecell}       % 美化单元格内容
\usepackage{tabularx}       % 自动调整列宽的表格
\usepackage{nicematrix}     % 增强矩阵排版
\usepackage{tikz}\usetikzlibrary{tikzmark, arrows.meta}   % 矢量图
\usepackage{float} 
\usepackage{xcolor}           % 颜色控制
\usepackage{graphicx}         % 图片插入
\usepackage{subcaption}

\usepackage{subcaption}
\usepackage{hyperref} % 超链接
\usepackage{cleveref}
\hypersetup{
  bookmarksnumbered = true,
  bookmarksopen=false,
  pdfborder=0 0 0,         % make all links invisible, so the pdf looks good when printed
  pdffitwindow=true,      % window fit to page when opened
  pdfnewwindow=true, % links in new window
  colorlinks=true,           % false: boxed links; true: colored links
  linkcolor=blue,            % color of internal links
  citecolor=magenta,    % color of links to bibliography
  filecolor=magenta,     % color of file links
  urlcolor=cyan              % color of external links
}

\newcommand{\algcomment}[1]{%
  \item[]\hskip-\algorithmicindent
  {\small\itshape // #1}%
}

\ifpdf
  \DeclareGraphicsExtensions{.eps,.pdf,.png,.jpg}
\else
  \DeclareGraphicsExtensions{.eps}
\fi

\newsiamremark{hypothesis}{Hypothesis}
\crefname{hypothesis}{Hypothesis}{Hypotheses}
\newsiamthm{claim}{Claim}
\newsiamremark{fact}{Fact}
\crefname{fact}{Fact}{Facts}
\makeatletter
\newcommand\figcaption{\def\@captype{figure}\caption} % 用于非浮动环境的标题
\newcommand\tabcaption{\def\@captype{table}\caption}
\makeatother

\theoremstyle{definition}
\headers{SCRINN}{Y. Huang, W. Su, N. Yi and P. Yin}

\title{Sparsely connected rank-inspired neural network\thanks{Submitted to the editors DATE.
\funding{This research was partially supported by China's National Key R \& D Programs (2024YFA1012600), NSFC Project (12431014).}}}

\author{
  Yunqing Huang\thanks{National Center for Applied Mathematics in Hunan, Key Laboratory of Intelligent Computing \& Information Processing of Ministry of Education, Xiangtan University, Xiangtan 411105, Hunan, P.R.China (\email{huangyq@xtu.edu.cn}).}
  \and
  Wu Su\thanks{Corresponding author. School of Mathematics and Computational Science, Xiangtan University, Xiangtan 411105, P.R.China (\email{202431510144@smail.xtu.edu.cn}).}
  \and
  Nianyu Yi\thanks{Hunan Key Laboratory for Computation and Simulation in Science and Engineering, Xiangtan University, Xiangtan 411105, P.R.China (\email{yinianyu@xtu.edu.cn}).}
  \and
  Peimeng Yin\thanks{Department of Mathematical Sciences, University of Texas at El Paso, El Paso, Texas 79968, USA (\email{pyin@utep.edu}).}
}
\usepackage{amsopn}